\documentclass[11pt]{article}

\title{Accurate and efficient splitting methods for dissipative particle dynamics}

\author{Xiaocheng Shang\footnote{Corresponding author. Email: \href{mailto:x.shang.1@bham.ac.uk}{x.shang.1@bham.ac.uk} } \\
\small{School of Mathematics, University of Birmingham, Edgbaston, Birmingham, B15 2TT, United Kingdom} }

\date{\today}

\usepackage{cite}

\usepackage{hyperref}
\usepackage{graphicx}

\usepackage{caption}
\DeclareGraphicsExtensions{.eps,.mps,.pdf,.jpg,.png}
\graphicspath{{figures/}{../figures/}}

\usepackage{fancyhdr}

\usepackage{amsmath}
\usepackage{bm}

\usepackage{amsfonts}

\usepackage{dsfont}

\usepackage[titletoc,toc,title]{appendix}

\usepackage[margin=3.0cm]{geometry}

\newcommand{\dd}{{\rm d}}
\newcommand{\kB}{k_{\mathrm{B}}}

\renewcommand{\vec}[1]{{\mathbf #1}}

\newcommand{\q}{{\vec{q}}}
\newcommand{\p}{{\vec{p}}}

\usepackage{mathrsfs}

\usepackage{soul}
\usepackage{color}

\usepackage{algorithm}
\usepackage{algorithmic}

\usepackage{url}
\usepackage{subfigure}

\usepackage{epstopdf}
\DeclareGraphicsExtensions{.eps,.mps,.pdf,.jpg,.png}

\usepackage{latexsym}
\usepackage{mathrsfs}
\usepackage{amssymb}
\usepackage{enumitem}

\usepackage{url}

\usepackage{array}
\newcolumntype{C}[1]{>{\centering\let\newline\\\arraybackslash\hspace{0pt}}m{#1}}

\begin{document}

\maketitle

\begin{abstract}
  We study numerical methods for dissipative particle dynamics (DPD), which is a system of stochastic differential equations and a popular stochastic momentum-conserving thermostat for simulating complex hydrodynamic behavior at mesoscales. We propose a new splitting method that is able to substantially improve the accuracy and efficiency of DPD simulations in a wide range of the friction coefficients, particularly in the extremely large friction limit that corresponds to a fluid-like Schmidt number, a key issue in DPD. Various numerical experiments on both equilibrium and transport properties are performed to demonstrate the superiority of the newly proposed method over popular alternative schemes in the literature.
\end{abstract}

\pagenumbering{arabic}

\section{Introduction}
\label{sec:Introduction}

Since its introduction by Hoogerbrugge and Koelman~\cite{Hoogerbrugge1992} at the end of the last century and due to its algorithmic simplicity and modelling versatility, DPD has become a popular stochastic momentum-conserving thermostat for simulating complex hydrodynamic behavior at a mesoscopic level~\cite{Groot1997,Soddemann2003}. Unlike individual molecules, DPD particles represent groups of fluid molecules and interact at short range with a soft potential. The coarse-grained descriptions allow the use of time and length scales that would otherwise be inaccessible by conventional molecular dynamics at microscales~\cite{Allen2017,Frenkel2001}. Therefore, DPD has been widely used in a large number of complex fluids and soft matter applications, including colloids~\cite{Koelman1993}, blood~\cite{Fedosov2013}, and polymers~\cite{Spenley2000} (see more applications in an excellent recent review~\cite{Espanol2017} and references therein).

\subsection{Formulation of DPD}
\label{subsec:Formulation}

Originally updated in discrete time steps, DPD was later reformulated by Espa\~{n}ol and Warren~\cite{Espanol1995} as a proper statistical mechanics model that is a system of It\={o} stochastic differential equations (SDEs). Consider an $N$-particle system evolving in dimension $d$ with position $\q_{i} \in \mathbb{R}^{d}$, momentum $\p_{i} \in \mathbb{R}^{d}$, and mass $m_{i} \in \mathbb{R}$ for $i=1,\dots,N$, the equations of motion for DPD particles are given by
\begin{subequations}\label{eq:DPD}
\begin{align}
  \dd \mathbf{q}_{i} &= m_{i}^{-1}\mathbf{p}_{i} \dd t \, , \label{eq:DPD_q} \\
  \dd \mathbf{p}_{i} &= \sum_{j\neq i} \mathbf{F}_{ij}^{\mathrm{C}}(r_{ij}) \dd t - \gamma \sum_{j\neq i} \omega^{\mathrm{D}}(r_{ij})(\mathbf{e}_{ij}\cdot \mathbf{v}_{ij})\mathbf{e}_{ij} \dd t + \sigma \sum_{j\neq i} \omega^{\mathrm{R}}(r_{ij})\mathbf{e}_{ij} \dd \mathrm{W}_{ij} \, , \label{eq:DPD_p}
\end{align}
\end{subequations}
where $\mathbf{F}_{ij}^{\mathrm{C}}$ is the conservative force acting along the line of centres, typically chosen as~\cite{Groot1997}
\begin{equation}\label{eq:Fc}
  \mathbf{F}_{ij}^{\mathrm{C}}(r_{ij}) =
  \begin{cases}
  a_{ij} (1-r_{ij}/r_{\mathrm{c}}) \mathbf{e}_{ij} \, , & r_{ij} < r_{\mathrm{c}} \, ;\\
  \quad \quad \quad 0 \, , & r_{ij} \geq r_{\mathrm{c}} \, ,
  \end{cases}
\end{equation}
which is the derivative of the soft pair potential energy,
\begin{equation}\label{eq:Pair_Potential}
  \varphi(r_{ij})=
  \begin{cases}
    a_{ij}r_{\mathrm{c}} \left(1-r_{ij}/r_{\mathrm{c}}\right)^{2}/2 \, , & r_{ij} < r_{\mathrm{c}} \, ;\\
    \quad \quad \quad \quad 0 \, , & r_{ij} \geq r_{\mathrm{c}} \, ,
  \end{cases}
\end{equation}
where parameter $a_{ij}$ represents the maximum repulsion strength between particles $i$ and $j$, $r_{ij} = |\mathbf{q}_{ij}| = |\mathbf{q}_{i}-\mathbf{q}_{j}|$ is the distance, and $r_{\mathrm{c}}$ denotes a certain cutoff radius. Moreover, $\mathbf{e}_{ij}=\mathbf{q}_{ij}/r_{ij}$ is the unit vector pointing from particle $j$ to particle $i$, $\mathbf{v}_{ij} = \mathbf{p}_{i}/m_{i}-\mathbf{p}_{j}/m_{j}$ denotes the relative velocity, and $\dd \mathrm{W}_{ij}=\dd \mathrm{W}_{ji}$ represents independent increments of the Wiener process with mean zero and variance $\dd t$. Parameters $\gamma$ and $\sigma$, representing the dissipative and random strengths, respectively, and position-dependent weight functions $\omega^{\mathrm{D}}$ and $\omega^{\mathrm{R}}$ are required to satisfy the following fluctuation-dissipation relations for DPD:
\begin{equation}\label{eq:FDT}
  \sigma^{2}=2\gamma \kB T \, , \quad \omega^{\mathrm{D}}(r_{ij})=\left[\omega^{\mathrm{R}}(r_{ij})\right]^{2} \, ,
\end{equation}
where $\kB$ is the Boltzmann constant and $T$ is the equilibrium temperature. Note that one of the two weight functions can be arbitrarily chosen, for instance, a popular choice of $\omega^{\mathrm{R}}$ is
\begin{equation}\label{eq:omega_R}
  \omega^{\mathrm{R}}(r_{ij})=
  \begin{cases}
  1-r_{ij}/r_{\mathrm{c}} \, , & r_{ij}<r_{\mathrm{c}} \, ;\\
  \quad \quad 0 \, , & r_{ij}\geq r_{\mathrm{c}} \, ,
  \end{cases}
\end{equation}
which fixes the other weight function $\omega^{\mathrm{D}}$. It can then be easily shown that the canonical ensemble is preserved with an invariant measure defined by the density
\begin{equation}\label{eq:rho_beta}
  \rho_{\beta}(\mathbf{q},\mathbf{p}) = Z^{-1}\exp \left(-\beta H(\mathbf{q},\mathbf{p}) \right) \, ,
\end{equation}
where $\beta^{-1}=\kB T$, $Z$ is a suitable normalizing constant (i.e., the partition function), and $H$ denotes the Hamiltonian defined as
\begin{equation}
  H(\mathbf{q},\mathbf{p}) = \sum_{i} \frac{\mathbf{p}_{i} \cdot \mathbf{p}_{i}}{2m_{i}} + U(\mathbf{q}) \, ,
\end{equation}
where $U$ is the potential energy
\begin{equation}
  U(\mathbf{q}) = \sum_{i} \sum_{j>i}\varphi(r_{ij}) \, .
\end{equation}
Due to the pairwise (or symmetric) nature of the interactions between particles and also the dependence on relative velocities, both linear and angular momenta are conserved in DPD. Moreover, DPD is an isotropic Galilean-invariant thermostat that preserves hydrodynamics~\cite{Soddemann2003,Allen2007,Moeendarbary2009}. Note that if the linear momentum is conserved, the density~\eqref{eq:rho_beta} should be replaced by
\begin{equation}
  \rho_{\beta}(\mathbf{q},\mathbf{p}) = \frac{1}{Z} \exp(-\beta H(\mathbf{q},\mathbf{p})) \times \delta\left[ \sum_i p_{i,x}- \pi_x \right] \delta\left [\sum_i p_{i,y} - \pi_y \right] \delta\left [\sum_i p_{i,z} - \pi_z \right] \, ,
\end{equation}
where $\boldsymbol{\pi}=(\pi_x,\pi_y,\pi_z)$ is the linear momentum vector. A similar modification would also be needed for the angular momentum conservation. However, the angular momentum will not be conserved if either periodic boundary conditions or Lees--Edwards boundary conditions~\cite{Lees1972}
are applied. Note also that it is highly nontrivial to show the ergodicity of DPD, which has only been demonstrated in the case of high particle density in one dimension by Shardlow and Yan~\cite{Shardlow2006}.

\subsection{Schmidt number issue in DPD}
\label{subsec:Schmidt}

The Schmidt number, $\mathrm{Sc}$, is an important quantity that characterizes the dynamical behavior of fluids. It is the ratio of the kinematic viscosity $\nu$ (also called ``momentum diffusivity'') to the diffusion coefficient $D$, and in a typical fluid flow (e.g., water), momentum is expected to be transported more rapidly than mass, leading to a Schmidt number of $O(10^{3})$. However, it has been pointed out that the standard DPD formulation described in Section~\ref{subsec:Formulation} (also known as the ``model B'' in the language of~\cite{Nikunen2003}), with a standard set of parameters (e.g.,~\cite{Allen2006} with $\gamma=4.5$), produces a gas-like Schmidt number of $O(1)$, resulting in concerns on the separation of the timescale for the propagation of hydrodynamic interactions and that for diffusion~\cite{Groot1997}. More precisely, the Schmidt number associated with the standard DPD formulation is approximated as
\begin{equation}\label{eq:Schmidt}
  \mathrm{Sc} = \frac{\nu}{D} \approx \frac{1}{2} + \frac{\left( 2\pi \gamma \rho_{\rm d} r^{4}_{\mathrm{c}} \right)^{2}}{70875\kB T} \, ,
\end{equation}
where $\rho_{\rm d}$ is the particle density. Although we can easily see from the approximation above that the most efficient way to increase the Schmidt number is to extend the cutoff radius $r_{\mathrm{c}}$, it is also obvious that it could result in a substantial computational overhead~\cite{Fan2006}. As discussed in~\cite{Litvinov2010}, one could generate a larger viscosity by increasing the stiffness of the conservative force or the particle density. However, the represented length scale decreases in both approaches, thereby contradicting the intended coarse-graining property of the DPD method. Alternatively, it has been suggested by Groot and Warren in~\cite{Groot1997} that a larger value of the friction coefficient $\gamma$ may be used, although the authors also mentioned that the largest stepsize usable may have to be reduced in order to maintain the temperature control. To this end, the aim of the current article is to explore if it is possible to design novel numerical integrators that would allow the use of large stepsizes in potentially very large friction limit, producing a fluid-like Schmidt number, while maintaining good control of the temperature and other important physical quantities. To the best of our knowledge, we are not aware of such studies within the standard DPD formulation in the literature.

Although considerable effort has been devoted to developing accurate and efficient numerical methods where large stepsizes can be used while maintaining good equilibrium properties, the corresponding Schmidt numbers were often gas-like. Following early examinations on the performance of various DPD integrators~\cite{Besold2000,Vattulainen2002,Nikunen2003,Chaudhri2010}, a number of popular methods (including the Lowe--Andersen (LA) thermostat~\cite{Lowe1999}, and its variant, the Nos\'{e}--Hoover--Lowe--Andersen (NHLA) thermostat~\cite{Stoyanov2005}) have been systematically compared in two recent studies~\cite{Leimkuhler2015,Leimkuhler2016a}. In the current article, we demonstrate that the accuracy and efficiency of DPD simulations can be substantially improved in a wide range of the friction coefficients, especially in the extremely large friction limit (i.e, $\gamma=450$) that corresponds to a fluid-like Schmidt number (i.e., $\mathrm{Sc} \approx 1016$). It is worth mentioning that the Schmidt number can be varied by either modifying the weight function~\cite{Fan2006,Yaghoubi2015} or using alternative approaches (e.g., LA and NHLA). However, in the current article we restrict our attention to the standard DPD formulation as it is by far the most popular and studied approach.

\subsection{Outline of the article}

The rest of the article is organized as follows. In Section~\ref{sec:Numerical_Methods}, we first describe two popular integration methods for DPD, followed by the derivations of a new promising scheme. We also theoretically demonstrate the order of convergence for the newly proposed method that successively integrates the dissipative and random forces based on interacting pairs. A variety of numerical experiments are performed in Section~\ref{sec:Numerical_Experiments} to compare all the schemes described in the article with a wide range of the friction coefficients (we restrict our attention to the extremely large friction limit that corresponds to a fluid-like Schmidt number for transport properties, e.g., shear viscosity, in Section~\ref{subsubsec:Shear_Viscosity}). Our findings are summarized in Section~\ref{sec:Conclusions}.

\section{Numerical methods for DPD}
\label{sec:Numerical_Methods}

In discussing the numerical integration of DPD, it is more convenient to rewrite the DPD system~\eqref{eq:DPD} to a more compact form:
\begin{subequations}\label{eq:DPD_compact}
\begin{align}
  \dd \mathbf{q}_{i} &= m_{i}^{-1}\mathbf{p}_{i} \dd t \, , \\
  \dd \mathbf{p}_{i} &= \mathbf{F}^{\mathrm{C}}_{i} \dd t +  \mathbf{F}^{\mathrm{D}}_{i} \dd t + \dd \mathbf{F}^{\mathrm{R}}_{i} \, ,
\end{align}
\end{subequations}
where $\mathbf{F}^{\mathrm{C}}_{i}$, $\mathbf{F}^{\mathrm{D}}_{i}$, and $\dd \mathbf{F}^{\mathrm{R}}_{i}$ respectively represent the total conservative, dissipative, and random forces acting on particle $i$,
\begin{subequations}\label{eq:F}
\begin{align}
  \mathbf{F}^{\mathrm{C}}_{i}(\mathbf{q}) &=
  \sum_{j\neq i} \mathbf{F}_{ij}^{\mathrm{C}}(r_{ij}) = -\nabla_{\mathbf{q}_{i}}U(\mathbf{q})
  \, , \label{eq:F_C} \\
  \mathbf{F}^{\mathrm{D}}_{i}(\mathbf{q}, \mathbf{p}) &= - \gamma \sum_{j\neq i} \omega^{\mathrm{D}}(r_{ij})(\mathbf{e}_{ij}\cdot \mathbf{v}_{ij})\mathbf{e}_{ij} \, , \label{eq:F_D} \\
  \dd \mathbf{F}^{\mathrm{R}}_{i}(\mathbf{q}) &= \sigma \sum_{j\neq i} \omega^{\mathrm{R}}(r_{ij})\mathbf{e}_{ij} \dd \mathrm{W}_{ij} \, . \label{eq:F_R}
\end{align}
\end{subequations}

\subsection{The velocity Verlet method}

Despite the significant advancement of numerical integrators for DPD in the past two decades (see a comprehensive review of a large number of numerical methods and comparisons of several popular schemes in~\cite{Leimkuhler2015}), the velocity Verlet (VV) method~\cite{Besold2000} remains one of the most popular DPD integrators in popular software packages (e.g., LAMMPS~\cite{Plimpton1995} and DL\_MESO~\cite{Seaton2013}) due to its ease of implementation (particularly in parallel computing). The integration steps of the VV method read
\begin{subequations}\label{eq:VV}
\begin{align}
  \mathbf{p}_{i}^{n+1/2} &= \mathbf{p}_{i}^{n} + \left[ \Delta t \mathbf{F}^{\mathrm{C}}_{i}(\mathbf{q}^{n}) + \Delta t \mathbf{F}^{\mathrm{D}}_{i}(\mathbf{q}^{n},\mathbf{p}^{n}) + \dd \mathbf{F}^{\mathrm{R}}_{i}(\mathbf{q}^{n}) \right] /2 \, ,  \label{eq:VV_1} \\
  \mathbf{q}_{i}^{n+1} &=  \mathbf{q}_{i}^{n} + \Delta t m_{i}^{-1}\mathbf{p}_{i}^{n+1/2} \, , \label{eq:VV_2} \\
  \mathbf{p}_{i}^{n+1} &= \mathbf{p}_{i}^{n+1/2} + \left[ \Delta t \mathbf{F}^{\mathrm{C}}_{i}(\mathbf{q}^{n+1}) + \Delta t \mathbf{F}^{\mathrm{D}}_{i}(\mathbf{q}^{n+1},\mathbf{p}^{n+1/2}) + \dd \mathbf{F}^{\mathrm{R}}_{i}(\mathbf{q}^{n+1}) \right] /2 \, . \label{eq:VV_3}
\end{align}
\end{subequations}
where $\Delta t$ is the integration stepsize and $\dd \mathrm{W}_{ij}$ in the random force~\eqref{eq:F_R} is replaced by $\sqrt{\Delta t}\mathrm{R}_{ij}$ with $\mathrm{R}_{ij}$ being a normally distributed variable with zero mean and unit variance. Note that all the forces in~\eqref{eq:F} need to be computed only once in~\eqref{eq:VV_3} and are reused in the subsequent step, otherwise the factor associated with the random force (i.e., $\sqrt{\Delta t}/2$) would be different. The VV method is also known as MD-VV in~\cite{Vattulainen2002} and DL\_MESO, or the modified Verlet method in~\cite{Shardlow2003}, and is equivalent to the GW integrator of Groot and Warren~\cite{Groot1997} when the variable factor is fixed as one half (i.e., $\lambda=1/2$).

It is well known that the standard VV method in molecular dynamics is second order~\cite{Leimkuhler2005,Hairer2006}. However, due to the presence of the dissipative and random forces in DPD, only first order convergence to the invariant measure is expected in the DPD context~\cite{Shardlow2003}. It is worth mentioning a variant of the VV method, that is, the DPD-VV method~\cite{Besold2000}, also included in DL\_MESO. DPD-VV only differs from VV in performing an additional update of the dissipative forces at the end of the integration steps. Although it has been claimed in~\cite{Nikunen2003} that the additional update could improve the performance considerably (see also good overall performance of DPD-VV in~\cite{Besold2000,Vattulainen2002}), we have observed in~\cite{Leimkuhler2015} that the performance of DPD-VV is rather similar to that of the splitting method proposed by Shardlow~\cite{Shardlow2003}. Given that Shardlow's splitting method has been the most recommended DPD integrator in the literature~\cite{Nikunen2003,Lisal2011,Espanol2017}, it was often used to represent the standard DPD formulation in comparison studies~\cite{Leimkuhler2016a,Shang2017}. Therefore, we will include it in our numerical experiments in Section~\ref{sec:Numerical_Experiments}, while excluding DPD-VV.

\subsection{Shardlow's splitting method: DPD-S1}
\label{subsec:Shardlow}

Splitting methods have been widely used in a range of systems, including Hamiltonian dynamics~\cite{Leimkuhler2005,Hairer2006}, dissipative systems~\cite{Shang2018}, and stochastic dynamics~\cite{Leimkuhler2015a,Shang2017}. It was Shardlow who first adopted and systematically examined the techniques in the DPD context. More specifically, the vector field of the DPD system~\eqref{eq:DPD_compact} is decomposed into three parts, which we label as A, B, and O:
\begin{equation}
  \dd \left[ \begin{array}{c} \mathbf{q}_{i} \\ \mathbf{p}_{i} \end{array} \right] = \underbrace{\left[ \begin{array}{c} m_{i}^{-1}\mathbf{p}_{i} \\ \mathbf{0} \end{array} \right] \dd t}_\mathrm{A} + \underbrace{\left[ \begin{array}{c} \mathbf{0} \\ \mathbf{F}^{\mathrm{C}}_{i} \end{array} \right] \dd t}_\mathrm{B} + \underbrace{\left[ \begin{array}{c} \mathbf{0} \\ \mathbf{F}^{\mathrm{D}}_{i} \dd t + \dd \mathbf{F}^{\mathrm{R}}_{i} \end{array} \right]}_\mathrm{O} \, ,
\end{equation}
where the first two splitting pieces (A and B) represent the Hamiltonian (or deterministic) part of the system and each of the pieces can be solved exactly, while the remaining O part, with positions fixed, is an Ornstein--Uhlenbeck (OU) process. In describing splitting methods, we use the formal notation of the generator of the diffusion as in~\cite{DeFabritiis2006,Serrano2006,Thalmann2007}. The generators for each part of the system may be written out as follows:
\begin{subequations}
\begin{align}
  \mathcal{L}_\mathrm{A} &= \sum_{i}\frac{\mathbf{p}_{i}}{m_{i}} \cdot \nabla_{\mathbf{q}_{i}} \, , \label{eq:generator_A} \\
  \mathcal{L}_\mathrm{B} &= \sum_{i} \mathbf{F}_{i}^{\mathrm{C}} \cdot \nabla_{\mathbf{p}_{i}}
  = -\sum_{i} \nabla_{\mathbf{q}_{i}}U(\mathbf{q}) \cdot \nabla_{\mathbf{p}_{i}}
  \, , \label{eq:generator_B} \\
  \mathcal{L}_\mathrm{O} &= \sum_{i} \sum_{j>i}  \mathcal{L}_{\mathrm{O}_{i,j}} \, , \label{eq:generator_O}
\end{align}
\end{subequations}
where
\begin{equation}\label{eq:generator_O_ij}
  \mathcal{L}_{\mathrm{O}_{i,j}} = \left[ - \gamma \omega^{\mathrm{D}}(r_{ij})(\mathbf{e}_{ij} \cdot \mathbf{v}_{ij}) + \frac{\sigma^{2}}{2} \left[ \omega^{\mathrm{R}}(r_{ij}) \right]^{2} \mathbf{e}_{ij} \cdot \left( \nabla_{\mathbf{p}_{i}} - \nabla_{\mathbf{p}_{j}} \right) \right] \mathbf{e}_{ij} \cdot \left( \nabla_{\mathbf{p}_{i}} - \nabla_{\mathbf{p}_{j}} \right) \, .
\end{equation}
The generator for the DPD system thus can be written as
\begin{equation}\label{eq:generator_DPD}
  \mathcal{L}_\mathrm{DPD} = \mathcal{L}_\mathrm{A} + \mathcal{L}_\mathrm{B} + \mathcal{L}_\mathrm{O} \, .
\end{equation}
Moreover, the flow map (or phase space propagator) of the system may be given by the shorthand notation
\begin{equation}
  \mathcal{F}_{t} = e^{t \mathcal{L}_{\mathrm{DPD}}} \, ,
\end{equation}
where the exponential map is used to formally denote the solution operator. Furthermore, approximations of $\mathcal{F}_{t}$ may be obtained as products (taken in different arrangements) of exponentials of the various splitting terms. For instance, the phase space propagation of Shardlow's S1 splitting method~\cite{Shardlow2003}, termed DPD-S1, can be written as
\begin{equation}\label{eq:Propagator_DPD_S1}
  \exp\left(\Delta t \hat{\mathcal{L}}_\mathrm{DPD-S1} \right)
  = \exp\left(\Delta t \hat{\mathcal{L}}_\mathrm{O}\right) \exp\left(\frac{\Delta t}{2}\mathcal{L}_\mathrm{B}\right) \exp\left(\Delta t \mathcal{L}_\mathrm{A}\right) \exp\left(\frac{\Delta t}{2}\mathcal{L}_\mathrm{B}\right) \, ,
\end{equation}
where $\exp\left(\Delta t\mathcal{L}_f\right)$ denotes the phase space propagator associated with the corresponding vector field $f$. Note that the steplengths associated with various operations are uniform and span the interval $\Delta t$. Therefore the B step in~\eqref{eq:Propagator_DPD_S1} is taken with a steplength of $\Delta t/2$, while a steplength of $\Delta t$ is associated with either the A or O step. It is worth pointing out that, in dealing with pairwise interactions in the OU process, it is desirable to further split the vector field O into each interacting pair. Therefore, the propagation of the O part in~\eqref{eq:Propagator_DPD_S1} should be more explicitly defined as
\begin{equation}\label{eq:Propagator_O}
  \exp\left(\Delta t \hat{\mathcal{L}}_\mathrm{O}\right) = \exp\left(\Delta t \hat{\mathcal{L}}_{\mathrm{O}_{N-1,N}}\right) \dots \exp\left(\Delta t \hat{\mathcal{L}}_{\mathrm{O}_{1,3}}\right) \exp\left(\Delta t \hat{\mathcal{L}}_{\mathrm{O}_{1,2}}\right) \, .
\end{equation}
Since the method of Br\"{u}nger, Brooks, and Karplus (BBK)~\cite{Brunger1984} used for successively integrating each interacting pair in the O part has been shown to produce weak second order \mbox{approximations}~\cite{Shardlow2003}, the propagation associated with each interacting pair may be given by
\begin{equation}\label{eq:Propagator_O_ij}
  \exp\left(\Delta t \hat{\mathcal{L}}_{\mathrm{O}_{i,j}}\right) = \exp\left(\Delta t \left[ \mathcal{L}_{\mathrm{O}_{i,j}} + O(\Delta t^{2}) \right] \right) \, .
\end{equation}
where $\mathcal{L}_{\mathrm{O}_{i,j}}$ is defined in~\eqref{eq:generator_O_ij}. Overall, the BBK method is successively used for each interacting pair in the OU process (part O), followed by the velocity Verlet method~\cite{Verlet1967,Melchionna2007} for the Hamiltonian part where both A and B parts are solved exactly:
\\
\noindent \textbf{Step 1}: for each interacting pair within cutoff radius ($r_{ij}<r_{\mathrm{c}}$), in a successive manner~\cite{Larentzos2014},
\begin{subequations}
\begin{align}
  \mathbf{p}_{i}^{n+1/4} & = \mathbf{p}_{i}^{n} - K_{ij} (\mathbf{e}_{ij}^{n} \cdot \mathbf{v}_{ij}^{n}) \mathbf{e}_{ij}^{n} + {\bf J}_{ij} \, , \\
  \mathbf{p}_{j}^{n+1/4} & = \mathbf{p}_{j}^{n} + K_{ij} (\mathbf{e}_{ij}^{n} \cdot \mathbf{v}_{ij}^{n}) \mathbf{e}_{ij}^{n} - {\bf J}_{ij}\, , \\
  \mathbf{p}_{i}^{n+2/4} & = \mathbf{p}_{i}^{n+1/4} + {\bf J}_{ij} - \frac{K_{ij}}{1+2K_{ij}} \left[ (\mathbf{e}_{ij}^{n} \cdot \mathbf{v}_{ij}^{n+1/4}) \mathbf{e}_{ij}^{n} + 2{\bf J}_{ij} \right], \\
  \mathbf{p}_{j}^{n+2/4} & = \mathbf{p}_{j}^{n+1/4} - {\bf J}_{ij} + \frac{K_{ij}}{1+2K_{ij}} \left[ (\mathbf{e}_{ij}^{n} \cdot \mathbf{v}_{ij}^{n+1/4}) \mathbf{e}_{ij}^{n} + 2{\bf J}_{ij} \right],
\end{align}
\end{subequations}
where $K_{ij} = \gamma \omega^{\mathrm{D}}(r_{ij}^{n}) \Delta t/2 $ and ${\bf J}_{ij} = \sigma \omega^{\mathrm{R}}(r_{ij}^{n}) \mathbf{e}_{ij}^{n} \sqrt{\Delta t} \mathrm{R}_{ij}^{n}/2$.
\\
\noindent \textbf{Step 2}: for each particle $i$,
\begin{subequations}
\begin{align}
  \mathbf{p}_{i}^{n+3/4} & = \mathbf{p}_{i}^{n+2/4} + (\Delta t/2) \mathbf{F}_{i}^{\mathrm{C}}(\mathbf{q}^{n}) \, , \\
  \mathbf{q}_{i}^{n+1} & = \mathbf{q}_{i}^{n} + \Delta t m_{i}^{-1} \mathbf{p}_{i}^{n+3/4} \, , \\
  \mathbf{p}_{i}^{n+1} & = \mathbf{p}_{i}^{n+3/4} + (\Delta t/2) \mathbf{F}_{i}^{\mathrm{C}}(\mathbf{q}^{n+1}) \, . \label{eq:DPD_S1_S2-3}
\end{align}
\end{subequations}
Note that the conservative force $\mathbf{F}_{i}^{\mathrm{C}}$ needs to be computed only once in~\eqref{eq:DPD_S1_S2-3}, where the Verlet neighbor lists~\cite{Verlet1967} are also updated. The interacting pairs needed in the subsequent step can then be easily identified from the lists.

As a shorthand, we may term the DPD-S1 method~\eqref{eq:Propagator_DPD_S1} OBAB (similarly, the S2 method of Shardlow~\cite{Shardlow2003} would be equivalent to OBABO in the same language). It has been shown in~\cite{Shardlow2003,Chaudhri2010} that the accuracy of both methods are very close to each other in a number of physical quantities. Given that the S1 method is more efficient than S2~\cite{Shardlow2003}, in what follows only the S1 method will be examined as in~\cite{Nikunen2003,Thalmann2007,Lisal2011,Larentzos2014,Leimkuhler2015,Shang2017}.

\subsection{The ABOBA method}
\label{subsec:ABOBA}

Instead of using the BBK method as in DPD-S1, it may be more desirable to integrate each interacting pair in the OU process (part O) analytically as in the pairwise adaptive Langevin thermostat~\cite{Leimkuhler2016a}. It is worth mentioning that the analytical solution of the OU process is often preferred in Langevin dynamics~\cite{Chandrasekhar1943,Allen2017}. More precisely, for each interacting pair, $i$ and $j$ ($j>i$), subtracting $\dd \mathbf{v}_{j}$ from $\dd \mathbf{v}_{i}$ and multiplying the unit vector $\mathbf{e}_{ij}$ on both sides yields
\begin{equation}\label{eq:ABOBA_exact}
  m_{ij}\dd v_{ij} = - \gamma \omega^{\mathrm{D}}(r_{ij}) v_{ij} \dd t + \sigma \omega^{\mathrm{R}}(r_{ij}) \dd \mathrm{W}_{ij} \, ,
\end{equation}
where $m_{ij}=m_{i}m_{j}/(m_{i}+m_{j})$ is the ``reduced mass'' and $v_{ij}=\mathbf{e}_{ij} \cdot \mathbf{v}_{ij}$. The above equation is an OU process with the exact (in the sense of distributional fidelity) solution~\cite{Kloeden1992}
\begin{equation}\label{eq:ABOBA_exact_v}
  v_{ij}(t) = e^{-\tau t} v_{ij}(0) + \frac{\sigma \omega^{\mathrm{R}}}{m_{ij}} \sqrt{ \frac{1-e^{-2\tau t}}{2\tau} }\mathrm{R}_{ij} \, ,
\end{equation}
where $\tau=\gamma \omega^{\mathrm{D}} / m_{ij}$ and $v_{ij}(0)$ is the initial relative velocity. Thus the velocity increment can be obtained as
\begin{equation}
  \Delta v_{ij} = v_{ij}(t) - v_{ij}(0) = v_{ij}(0) \left( e^{-\tau t} - 1 \right) + \frac{\sigma \omega^{\mathrm{R}}}{m_{ij}} \sqrt{ \frac{1-e^{-2\tau t}}{2\tau} }\mathrm{R}_{ij} \, ,
\end{equation}
and the corresponding momenta can be updated as follows:
\begin{subequations}\label{eq:ABOBA_exact_momenta}
\begin{align}
  \mathbf{p}_{i} &\leftarrow \mathbf{p}_{i} + m_{ij} \Delta v_{ij} \mathbf{e}_{ij} \, , \\
  \mathbf{p}_{j} &\leftarrow \mathbf{p}_{j} - m_{ij} \Delta v_{ij} \mathbf{e}_{ij} \, ,
\end{align}
\end{subequations}
which defines the propagator, $\exp({\Delta t \mathcal{L}_{\mathrm{O}_{i,j}}})$, for each interacting pair.

Given that the BAOAB method is the best-performing method in terms of sampling configurational quantities in Langevin dynamics~\cite{Leimkuhler2013,Leimkuhler2013a,Leimkuhler2015b} (see also a recent comprehensive study on its time correlations~\cite{Shang2019}), it is worth mentioning that a BAOAB method in the DPD context can be easily constructed following the procedures in Section~\ref{subsec:ABOBA}. However, it is important to note that in the DPD context, when updating the OU process in the BAOAB method, all the distances between (updated) interacting pairs have to be recalculated, the cost of which is essentially the same as another force calculation. In what follows we propose a new method, which we term the ABOBA method, where the recalculation can be easily avoided with the help of the Verlet neighbor lists. The propagator of the ABOBA method can be written as
\begin{equation}\label{eq:Propagator_ABOBA}
  \exp\left(\Delta t \hat{\mathcal{L}}_\mathrm{ABOBA}\right) = \exp\left(\frac{\Delta t}{2}\mathcal{L}_\mathrm{A}\right) \exp\left(\frac{\Delta t}{2}\mathcal{L}_\mathrm{B}\right) \exp\left(\Delta t \hat{\mathcal{L}}_\mathrm{O}\right)
  \exp\left(\frac{\Delta t}{2}\mathcal{L}_\mathrm{B}\right)
  \exp\left(\frac{\Delta t}{2}\mathcal{L}_\mathrm{A}\right) \, ,
\end{equation}
where the propagation of the O part is similarly defined as in~\eqref{eq:Propagator_O} except each interacting pair is solved exactly as demonstrated in this section. Note that one may wish to reverse the order of the interacting pairs in the O part, which will affect neither its overall performance nor the order of convergence to the invariant measure (see more discussions in Section~\ref{subsec:Accuracy}). The detailed integration steps of the ABOBA method read:
\\
\noindent \textbf{Step 1}: for each particle $i$,
\begin{align}
  \mathbf{q}_{i}^{n+1/2} & = \mathbf{q}_{i}^{n} + (\Delta t/2) m_{i}^{-1} \mathbf{p}_{i}^{n} \, , \\
  \mathbf{p}_{i}^{n+1/3} & = \mathbf{p}_{i}^{n} + (\Delta t/2) \mathbf{F}_{i}^{\mathrm{C}}(\mathbf{q}^{n+1/2}) \, , \label{eq:ABOBA_S1-2}
\end{align}
\\
\noindent \textbf{Step 2}: for each interacting pair within cutoff radius ($r_{ij}<r_{\mathrm{c}}$), in a successive manner,
\begin{subequations}
\begin{align}
  \mathbf{p}^{n+2/3}_{i} &= \mathbf{p}^{n+1/3}_{i} + m_{ij} \Delta v_{ij} \mathbf{e}^{n+1/2}_{ij} \, , \\
  \mathbf{p}^{n+2/3}_{j} &= \mathbf{p}^{n+1/3}_{j} - m_{ij} \Delta v_{ij} \mathbf{e}^{n+1/2}_{ij} \, ,
\end{align}
\end{subequations}
with
\begin{equation}
  \Delta v_{ij} = \left[ \mathbf{e}^{n+1/2}_{ij} \cdot \mathbf{v}^{n+1/3}_{ij} \right] \left( e^{-\tau \Delta t} - 1 \right) + \frac{\sigma \omega^{\mathrm{R}}(r^{n+1/2}_{ij})}{m_{ij}} \sqrt{ \frac{1-e^{-2\tau \Delta t}}{2\tau} }\mathrm{R}^{n}_{ij} \, ,
\end{equation}
where $\tau = \gamma \omega^{\mathrm{D}}(r^{n+1/2}_{ij}) / m_{ij}$.
\\
\noindent \textbf{Step 3}: for each particle $i$,
\begin{align}
  \mathbf{p}_{i}^{n+1} & = \mathbf{p}_{i}^{n+2/3} + (\Delta t/2) \mathbf{F}_{i}^{\mathrm{C}}(\mathbf{q}^{n+1/2}) \, , \\
  \mathbf{q}_{i}^{n+1} & = \mathbf{q}_{i}^{n+1/2} + (\Delta t/2) m_{i}^{-1} \mathbf{p}_{i}^{n+1} \, .
\end{align}
Note that the conservative force $\mathbf{F}_{i}^{\mathrm{C}}$ needs to be computed only once in~\eqref{eq:ABOBA_S1-2}, where the Verlet neighbor lists are also updated. The interacting pairs needed in the subsequent step can then be easily identified from the lists. Note also that if we switch off the O part, the ABOBA method reduces to the position Verlet method~\cite{Melchionna2007} for the Hamiltonian part where again A and B parts are solved exactly.

Denoting the operator of part S as the sum of the operators of B and O parts defined in Section~\ref{subsec:Shardlow} (i.e., $\mathcal{L}_\mathrm{S} = \mathcal{L}_\mathrm{B} + \mathcal{L}_\mathrm{O}$), it is also possible to integrate each interacting pair in the S part analytically as in the DPD-Trotter method~\cite{DeFabritiis2006,Serrano2006}. In the same language as in ABOBA, the DPD-Trotter method would be equivalent to SAS where each interacting pair in the S part is solved exactly in a successive manner. Moreover, a similar method of ASA could be easily constructed. While, with a friction coefficient of $\gamma=4.5$, the performance of DPD-Trotter/SAS is very similar to that of DPD-S1 as reported in~\cite{DeFabritiis2006,Leimkuhler2015}, the ASA method performs much better. However, as we increase the value of the friction coefficient, the performance of both SAS and ASA methods deteriorate significantly as the influence of the conservative force becomes less and less in both methods---in the limit of $\gamma \rightarrow \infty$, both methods effectively correspond to the integrations of an ideal gas (sometimes termed ``ideal DPD fluid'' within the DPD framework~\cite{Vattulainen2002}). However, the conservative force plays a crucial role in guiding the movements of the particles, and thus should not be ``neglected''. Moreover, since we are interested in the extremely large friction limit in order to achieve a fluid-like Schmidt number (see more discussions in Section~\ref{subsec:Schmidt}), we will include neither of the methods for comparisons in the current article.

\subsection{Accuracy of equilibrium averages}
\label{subsec:Accuracy}


The framework of long-time Talay--Tubaro expansion has been widely used in the analysis of the accuracy of ergodic averages (with respect to the invariant measure) in stochastic numerical methods~\cite{Talay1990,Debussche2012,Leimkuhler2013,Leimkuhler2013a,Abdulle2014a,Abdulle2014,Leimkuhler2013c,Leimkuhler2015a,Leimkuhler2015b}. While the DPD-S1 method~\eqref{eq:Propagator_DPD_S1} has been shown in~\cite{Leimkuhler2015} to have second order convergence to its invariant measure, in what follows we adopt the procedures to examine the order of convergence to the invariant measure of the newly proposed ABOBA method described in the previous section.

For a given splitting method described by $\mathcal{L}= \mathcal{L}_{\alpha} + \mathcal{L}_{\beta} + \dots + \mathcal{L}_{\zeta}$, we define its associated effective operator $\hat{\mathcal{L}}^{\dag}$ with stepsize $\Delta t$ as
\begin{equation}
  \exp\left(\Delta t \hat{\mathcal{L}}^{\dag}\right) = \exp\left(\Delta t \mathcal{L}_{\alpha}^{\dag}\right) \exp\left(\Delta t \mathcal{L}_{\beta}^{\dag}\right) \dots \exp\left(\Delta t \mathcal{L}_{\zeta}^{\dag}\right) \, ,
\end{equation}
where $\mathcal{L}_{\alpha}^{\dag}, \mathcal{L}_{\beta}^{\dag},\dots,\mathcal{L}_{\zeta}^{\dag}$ represent the corresponding Fokker--Planck operator associated with each subsystem. The effective operator $\hat{\mathcal{L}}^{\dag}$ can be viewed as a perturbation of the exact Fokker--Planck operator $\mathcal{L}^{\dag} = \mathcal{L}_{\alpha}^{\dag} + \mathcal{L}_{\beta}^{\dag} + \dots + \mathcal{L}_{\zeta}^{\dag}$ for the whole system:
\begin{equation}\label{eq:hat_L}
  \hat{\mathcal{L}}^{\dag} = \mathcal{L}^{\dag} + \Delta t \mathcal{L}^{\dag}_{1} + \Delta t^{2}\mathcal{L}^{\dag}_{2} + O(\Delta t^{3}) \, ,
\end{equation}
where perturbation operators $\mathcal{L}^{\dag}_{1}, \mathcal{L}^{\dag}_{2},\dots$ can be computed by using the Baker--Campbell--Hausdorff expansion~\cite{Hairer2006,Leimkuhler2005}. We also define the perturbed invariant measure $\hat{\rho}$ associated with the numerical method as an approximation of the target invariant measure $\rho_{\beta}$:
\begin{equation}\label{eq:hat_rho}
  \hat{\rho} = \rho_{\beta}\left[ 1 + \Delta t f_{1} + \Delta t^{2}f_{2} + \Delta t^{3}f_{3} + O(\Delta t^{4}) \right] \, ,
\end{equation}
where $f_{1}, f_{2},\dots$ are some correction functions. The average of each of those functions with respect to the target invariant measure is zero, i.e., $\langle f_{i} \rangle=0$. Subsequently, substituting $\hat{\mathcal{L}}^{\dag}$ and $\hat{\rho}$ into the stationary Fokker--Planck equation
\begin{equation}
  \hat{\mathcal{L}}^{\dag}\hat{\rho} = 0
\end{equation}
yields
\begin{equation}
  \left[ \mathcal{L}^{\dag} + \Delta t \mathcal{L}^{\dag}_{1} + \Delta t^{2}\mathcal{L}^{\dag}_{2} + O(\Delta t^{3}) \right] \left( \rho_{\beta} \left[ 1 + \Delta t f_{1} + \Delta t^{2} f_{2} + \Delta t^{3} f_{3} + O(\Delta t^{4}) \right] \right)=0 \, .
\end{equation}
Since the exact Fokker--Planck operator preserves the target invariant measure, i.e., $\mathcal{L}^{\dag}\rho_{\beta}=0$, we obtain
\begin{equation}\label{eq:error_analysis_PDE}
  \mathcal{L}^{\dag}(\rho_{\beta}f_{1}) = - \mathcal{L}^{\dag}_{1}\rho_{\beta}
\end{equation}
by equating first order terms in $\Delta t$. Although we are able to compute the perturbation operator $\mathcal{L}^{\dag}_{1}$ by using the Baker--Campbell--Hausdorff expansion for any particular splitting method, and subsequently its action on $\rho_{\beta}$, it is generally very hard to solve the above partial differential equation~\eqref{eq:error_analysis_PDE} in order to obtain the leading correction function $f_{1}$ in closed form (see examples in Langevin dynamics~\cite{Leimkuhler2013} and adaptive Langevin dynamics~\cite{Leimkuhler2015a}).

It is more convenient to work with the adjoint of the perturbed
generator associated with a particular splitting method.
For instance, in the ABOBA method~\eqref{eq:Propagator_ABOBA}, we have
\begin{equation}
  \exp\left(\Delta t \hat{\mathcal{L}}^{\dag}_\mathrm{ABOBA} \right)
  = \exp\left(\frac{\Delta t}{2}\mathcal{L}^{\dag}_\mathrm{A}\right) \exp\left(\frac{\Delta t}{2}\mathcal{L}^{\dag}_\mathrm{B}\right) \exp\left(\Delta t \hat{\mathcal{L}}^{\dag}_\mathrm{O}\right) \exp\left(\frac{\Delta t}{2}\mathcal{L}^{\dag}_\mathrm{B}\right) \exp\left(\frac{\Delta t}{2}\mathcal{L}^{\dag}_\mathrm{A}\right) \, ,
\end{equation}
where
\begin{equation}
  \exp\left(\Delta t \hat{\mathcal{L}}^{\dag}_\mathrm{O}\right) = \exp\left(\Delta t \mathcal{L}^{\dag}_{\mathrm{O}_{1,2}}\right) \exp\left(\Delta t \mathcal{L}^{\dag}_{\mathrm{O}_{1,3}}\right) \dots \exp\left(\Delta t \mathcal{L}^{\dag}_{\mathrm{O}_{N-1,N}}\right) \, .
\end{equation}
With the help of the Baker--Campbell--Hausdorff expansion, the effective operator associated with the overall method can be obtained as
\begin{equation}
  \hat{\mathcal{L}}^{\dag}_\mathrm{ABOBA} = \mathcal{L}^{\dag}_{\mathrm{A}} + \mathcal{L}^{\dag}_{\mathrm{B}} + \mathcal{L}^{\dag}_{\mathrm{O}} + \Delta t \mathcal{L}^{\dag}_{1,\mathrm{ABOBA}} + O(\Delta t^{2}) \, ,
\end{equation}
where
\begin{equation}\label{eq:L1_ABOBA}
\begin{aligned}
  \mathcal{L}^{\dag}_{1,\mathrm{ABOBA}} =& \, \frac{1}{2} \left[ \mathcal{L}^{\dag}_{\mathrm{O}_{1,2}}, \mathcal{L}^{\dag}_{\mathrm{O}_{1,3}} \right] + \frac{1}{2} \left[ \mathcal{L}^{\dag}_{\mathrm{O}_{1,2}} + \mathcal{L}^{\dag}_{\mathrm{O}_{1,3}}, \mathcal{L}^{\dag}_{\mathrm{O}_{1,4}} \right] + \dots \\
  &  + \frac{1}{2} \left[ \mathcal{L}^{\dag}_{\mathrm{O}_{1,2}} + \mathcal{L}^{\dag}_{\mathrm{O}_{1,3}} + \dots + \mathcal{L}^{\dag}_{\mathrm{O}_{N-2,N}}, \mathcal{L}^{\dag}_{\mathrm{O}_{N-1,N}} \right] \, .
\end{aligned}
\end{equation}
It can be easily shown that
\begin{equation}
   \mathcal{L}^{\dag}_{\mathrm{O}_{i,j}} \rho_{\beta} = 0 \, ,
\end{equation}
and subsequently
\begin{equation}
  \mathcal{L}^{\dag}_{1,\mathrm{ABOBA}} \rho_{\beta} = 0 \, .
\end{equation}
Therefore, the associated leading correction function in~\eqref{eq:error_analysis_PDE} must be zero, i.e.,
\begin{equation}
  f_{1,\mathrm{ABOBA}} = 0 \, .
\end{equation}
Given that higher order perturbations in~\eqref{eq:hat_rho} are not equal to zero in general, we have shown that the ABOBA method has second order convergence to its invariant measure.

\section{Numerical experiments}
\label{sec:Numerical_Experiments}

In this section, a variety of numerical experiments are conducted to systematically compare the newly proposed ABOBA method with alternative popular methods described in Section~\ref{sec:Numerical_Methods}.

\subsection{Simulation details}

In our numerical experiments, we adopted a standard set of parameters commonly used in algorithms tests as in~\cite{Groot1997,Allen2006}. Specifically, particle mass $m_{i}$, cutoff radius $r_{\mathrm{c}}$, and $\kB T$ were set to be unity. Although a particle density of $\rho_{\rm d}=4$ was initially used in~\cite{Groot1997} and later adopted in a number of studies~\cite{Vattulainen2002,Nikunen2003,Shardlow2003}, a smaller value of $\rho_{\rm d}=3$ was later suggested for efficiency reasons and thus was used throughout the current article. Subsequently, a repulsion parameter of $a_{ij} = 75\kB T/\rho_{\rm d} = 25$ was determined in order to match the compressibility of water~\cite{Groot1997}. Although a friction coefficient of $\gamma=4.5$ was widely used in algorithms tests, the corresponding Schmidt number was only $\mathrm{Sc}\approx 0.6$ according to~\eqref{eq:Schmidt}---a gas-like Schmidt number (see more discussions in Section~\ref{subsec:Schmidt}). While larger value of $\gamma=40.5$, which corresponds to $\mathrm{Sc}\approx 8.7$, was examined in~\cite{Leimkuhler2015,Leimkuhler2016a}, even larger values of the friction coefficient are needed to obtain fluid-like Schmidt numbers. Therefore, we also included friction coefficients of $\gamma=200$ and $\gamma=450$, corresponding to Schmidt numbers of $\mathrm{Sc}\approx 201$ and $\mathrm{Sc}\approx 1016$, respectively.

Moreover, a system of $N=500$ identical particles was simulated in a cubic box with periodic boundary conditions~\cite{Allen2017,Frenkel2001}, unless otherwise stated. While the initial positions of the particles were independent and identically distributed (i.i.d.) with a uniform distribution over the box, the initial momenta were i.i.d.\ normal random variables with mean zero and variance $\kB T$. Verlet neighbor lists~\cite{Verlet1967} were used wherever possible in order to reduce the computational cost as discussed in Section~\ref{sec:Numerical_Methods}.

\subsubsection{Equilibrium properties}
\label{subsubsec:Equilibrium}

Following~\cite{Leimkuhler2015,Leimkuhler2016a}, we measured the ``numerical efficiency'', defined as the ratio of the ``critical stepsize'' and the CPU time per step, of each method and then scaled it to that of the benchmark VV method, unless otherwise stated. The CPU time (in milliseconds) for the main integration steps (without calculating any physical quantities) was the time taken (on a Lenovo ThinkStation P330 Tiny) for the integration of a single time step of $\Delta t=0.05$ (averaged over 10,000 consecutive time steps). Note that Verlet neighbor lists~\cite{Verlet1967} were again used wherever possible. As in~\cite{Leimkuhler2015,Leimkuhler2016a}, the critical stepsize was determined as the stepsize that approximately corresponds to a 10\% relative error in the computed configurational temperature~\cite{Rugh1997,Braga2005,Allen2006,Travis2008} (unless otherwise stated), an observable function that depends solely on the positions. Moreover, the average of the computed configurational temperature in the canonical ensemble is expected to be precisely the target temperature:
\begin{equation}\label{eq:Config_Temp}
  \kB T = \frac{  \left\langle \nabla_{i} U(\mathbf{q}) \cdot \nabla_{i} U(\mathbf{q}) \right\rangle }{ \left\langle \nabla^{2}_{i}U(\mathbf{q}) \right\rangle } \, ,
\end{equation}
where $\nabla_{i}U$ and $\nabla^{2}_{i}U$ respectively represent the gradient and Laplacian of the potential energy $U$ with respect to the position of particle $i$ (see more discussions on the configurational temperature in~\cite{Leimkuhler2015,Leimkuhler2016a}). It should be noted that since the canonical momentum distribution is always Gaussian (and thus trivial to sample), as in Langevin dynamics~\cite{Leimkuhler2015b} we are far more interested in sampling configurational quantities. Thus the configurational temperature~\eqref{eq:Config_Temp} was chosen over the kinetic temperature that depends solely on the momenta. Importantly, good control of the configurational temperature appears to imply good performance in other physical quantities tested in Section~\ref{sec:Numerical_Experiments} (see more discussions on the reasoning in~\cite{Leimkuhler2015}). Moreover, it has been recommended in~\cite{Allen2006} that the configurational temperature~\eqref{eq:Config_Temp}, as a verification of equilibrium, should be measured and reported in DPD simulations.

In addition, we calculated the radial distribution function (RDF)~\cite{Allen2017,Frenkel2001}, often denoted as $g(r)$, which is another important configurational quantity in simulations, characterizing the structure of the system.

\subsubsection{Transport properties}
\label{subsubsec:Transport}

\begin{figure}[tb]
\centering
\includegraphics[scale=0.4]{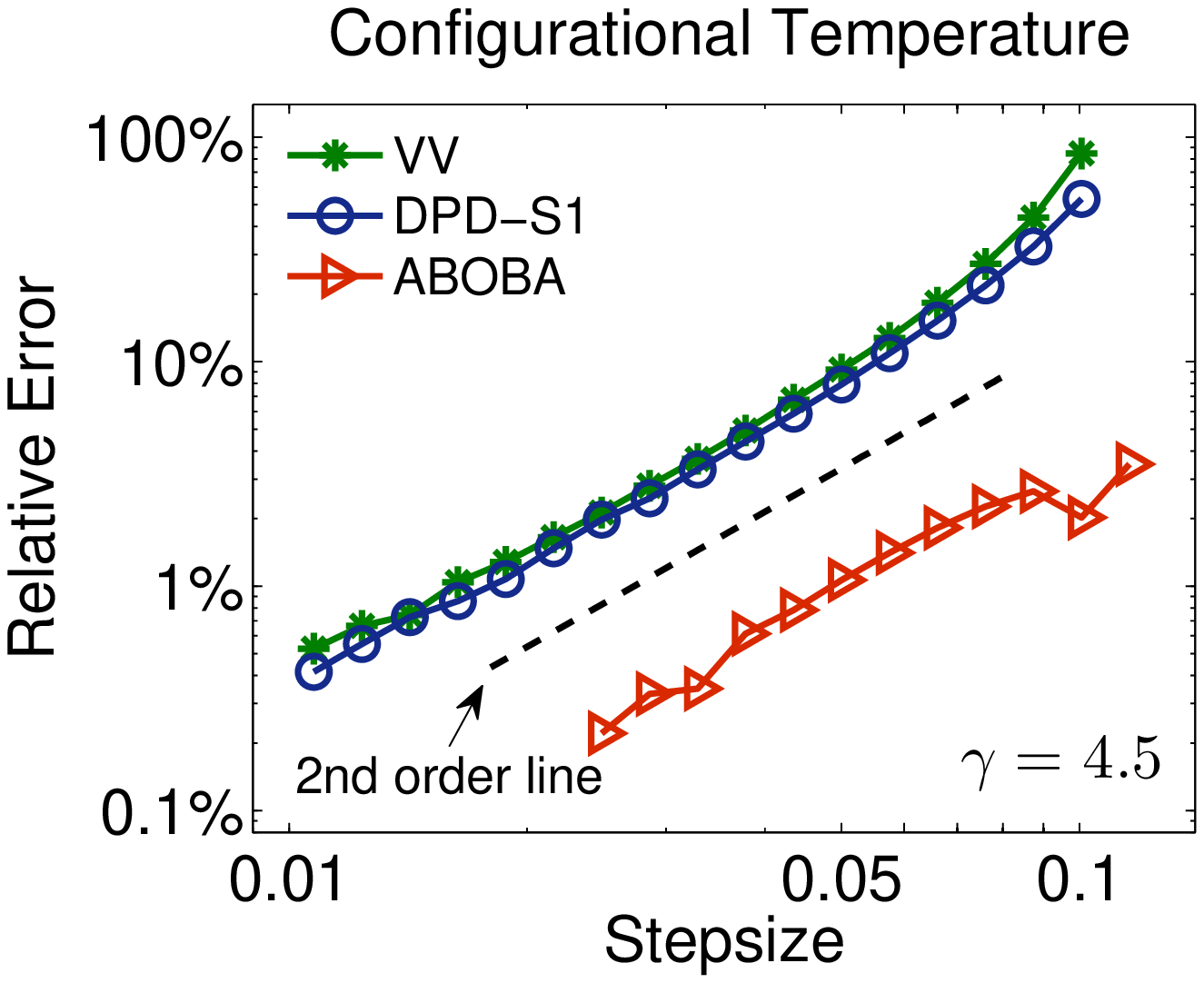}
\includegraphics[scale=0.4]{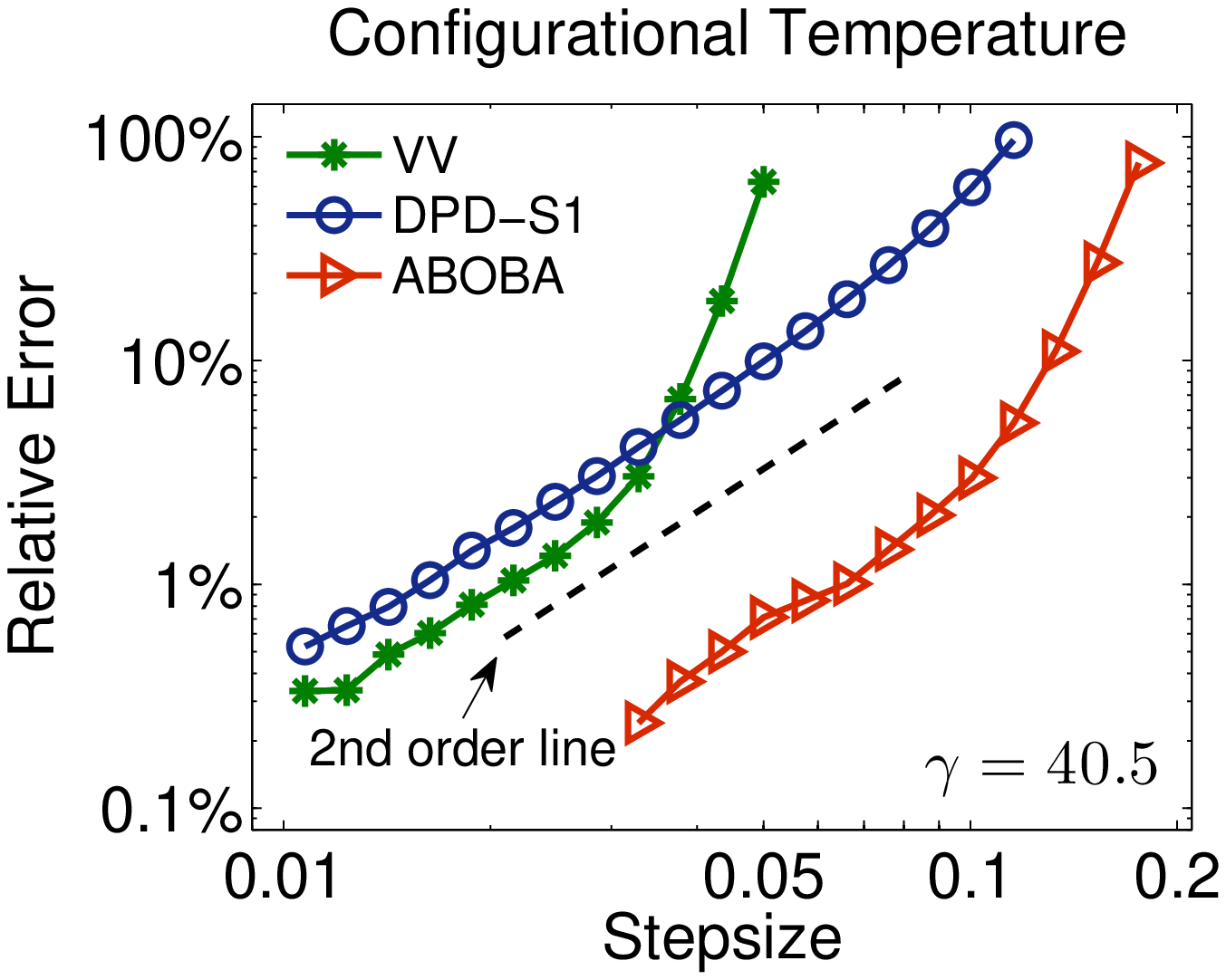}
\caption{\small Double logarithmic plot of the relative error in the computed configurational temperature~\eqref{eq:Config_Temp} against stepsize by using various numerical methods of the DPD system described in Section~\ref{sec:Numerical_Methods} with friction coefficients of $\gamma=4.5$ (left) and $\gamma=40.5$ (right). The system was simulated for 1000 reduced time units but only the last 80\% of the data were collected to calculate the static quantity in order to make sure the system was well equilibrated. Ten different runs were averaged to reduce the sampling errors. The stepsizes tested began at around $\Delta t=0.01$ and were increased incrementally by 15\% until all methods either started to show significant relative errors or became unstable (e.g., the VV method on the right panel at slightly over $\Delta t=0.05$). The dashed black line represents the second order convergence to the invariant measure.}
\label{fig:DPD_CT_Comp_g4d5+g40d5}
\end{figure}



A common and favourable approach to measure transport coefficients in particle-based methods is to employ planar Couette flow as a numerical ``viscometer''~\cite{Allen2017,Evans2008} (see a detailed discussion on extracting transport coefficients by a variety of approaches in~\cite{Shang2015a}). As a nonequilibrium method, a simple and steady shear flow is commonly generated in DPD via the well-known Lees--Edwards boundary conditions~\cite{Lees1972} in order to measure the shear viscosity numerically~\cite{Pagonabarraga1998,Backer2005,Chaudhri2010} (see also theoretical studies in~\cite{Marsh1997a,Groot1997,Evans1999,Masters1999}).
With Lees--Edwards boundary conditions, the primary cubic box remains centered at the origin as with normal periodic boundary conditions. However, a uniform shear velocity profile is generated with a streaming velocity that corresponds to the location of particle $i$~\cite{Evans2008}
\begin{equation}\label{eq:LEBC_Streaming_V}
  \mathbf{u}_i = \kappa({\bf q}_i\cdot{\bf e}^y)\mathbf{e}^{x} = \boldsymbol{\kappa}\cdot{\bf q}_i,
  \qquad \boldsymbol{\kappa} = \kappa\,{\bf e}^x\otimes{\bf e}^y
\end{equation}
where $\kappa$ is the shear rate defined as $\kappa \equiv \kappa_{xy} = \dd u^{x}/\dd y$, with $u^{x}$ being the macroscopic velocity in the $x$-direction, $\mathbf{e}^{x}$ and $\mathbf{e}^{y}$ respectively represent the unit vectors in the $x$- and $y$-direction, $\boldsymbol{\kappa}$ denotes the transposed velocity gradient tensor, and $\otimes$ is the dyadic product of two vectors. Note that Lees--Edwards boundary conditions are typically only applied in the $x$-direction, while the other directions ($y$ and $z$) remain with periodic boundary conditions. Particular care should be taken when implementing Lees--Edwards boundary conditions in pairwise thermostats, such as DPD, due to the position-dependence on both dissipative and random forces, this issue has been discussed in~\cite{Leimkuhler2016a} in order to maintain correct velocity profiles especially in large friction limits.

\begin{table}[t]
\centering
  \begin{tabular}{cccc}
  \hline
  \textbf{Method} & \textbf{Critical stepsize} & \textbf{CPU time} & \textbf{Scaled efficiency} \\ \hline
  VV     & 0.050 & 6.829 & 100.0\% \\
  DPD-S1 & 0.057 & 7.852 & 99.2\% \\
  ABOBA  & 0.116 & 7.742 & 204.6\% \\ \hline
  \end{tabular}
\caption[Table caption text]{\small Comparisons of the ``numerical efficiency'' of various numerical methods of the DPD system with a friction coefficient of $\gamma=4.5$. ``Critical stepsize'' is the stepsize beyond which the numerical method starts to show pronounced artifacts (i.e., 10\% relative error in the computed configurational temperature according to the left panel of Figure~\ref{fig:DPD_CT_Comp_g4d5+g40d5}). The numerical efficiency of each method was scaled to that of the benchmark VV method.}
\label{table:efficiency_g4d5}
\end{table}

\begin{table}[t]
\centering
  \begin{tabular}{cccc}
  \hline
  \textbf{Method} & \textbf{Critical stepsize} & \textbf{CPU time} & \textbf{Scaled efficiency} \\ \hline
  VV     & 0.038 & 6.829 & 100.0\% \\
  DPD-S1 & 0.050 & 7.852 & 114.4\% \\
  ABOBA  & 0.116 & 7.742 & 269.3\% \\ \hline
  \end{tabular}
\caption[Table caption text]{\small Comparisons of the ``numerical efficiency'' of various numerical methods of the DPD system with a large friction coefficient of $\gamma=40.5$, corresponding to the right panel of Figure~\ref{fig:DPD_CT_Comp_g4d5+g40d5}. The format of the table is the same as in Table~\ref{table:efficiency_g4d5}.}
\label{table:efficiency_g40d5}
\end{table}

The Irving--Kirkwood stress tensor~\cite{Irving1950} subject to Lees--Edwards boundary conditions can be written as
\begin{equation}\label{eq:Stree_Tensor}
    \boldsymbol{\sigma} = - \frac{1}{V} \left[ \sum_{i} m_{i}  \left( \mathbf{v}_{i}-\mathbf{u}_i \right) \otimes \left( \mathbf{v}_{i}-\mathbf{u}_i \right) + \sum_{i} \sum_{j>i} \mathbf{q}_{ij} \otimes \mathbf{F}_{ij} \right] \, ,
\end{equation}
where $V$ is the volume of the simulation box. While only the conservative force should be accounted for $\mathbf{F}_{ij}$ in Langevin dynamics since both the dissipative and random forces are averaged out, all three components of the force should be included for pairwise thermostats, including DPD. The shear viscosity can be extracted at finite rates as
\begin{equation}\label{eq:Shear_Viscosity}
    \eta = \frac{\langle \sigma_{xy} \rangle}{\kappa} \, ,
\end{equation}
where $\sigma_{xy}$ represents the shear stress and is the off-diagonal $xy$-component of the symmetric stress tensor $\boldsymbol{\sigma}$~\eqref{eq:Stree_Tensor}. It is worth mentioning that the zero shear viscosity $\eta_0=\lim_{\kappa\rightarrow 0}\eta$ can be obtained by extrapolation from~\eqref{eq:Shear_Viscosity}. Note also that $\eta_0$ can be alternatively calculated by integrating the stress-stress autocorrelation function (i.e., the Green--Kubo formulas~\cite{Green1954,Kubo1957}) as in~\cite{Chaudhri2010}. However, it is well documented that those equilibrium approaches are subject to significant statistical error and thus not preferred in practice~\cite{Allen2017,Shang2015a}.

The control of another important transport coefficient, the diffusion coefficient $D$, which is proportional to the integral of the unnormalized velocity autocorrelation function, has also been investigated. It appears that the accuracy of the approximation of the diffusion coefficient is highly sensitive to the stepsizes used, particularly in the extremely large friction limit. That is, a very small stepsize has to be used in order to achieve a reasonably accurate approximation of the diffusion coefficient. Since we are more interested in using large stepsizes, the results of the diffusion coefficient will not be presented.

\subsection{Numerical results}
\label{subsec:Numerical_Results}

\subsubsection{Configurational temperature}
\label{subsubsec:CT}

\begin{figure}[t]
\centering
\includegraphics[scale=0.4]{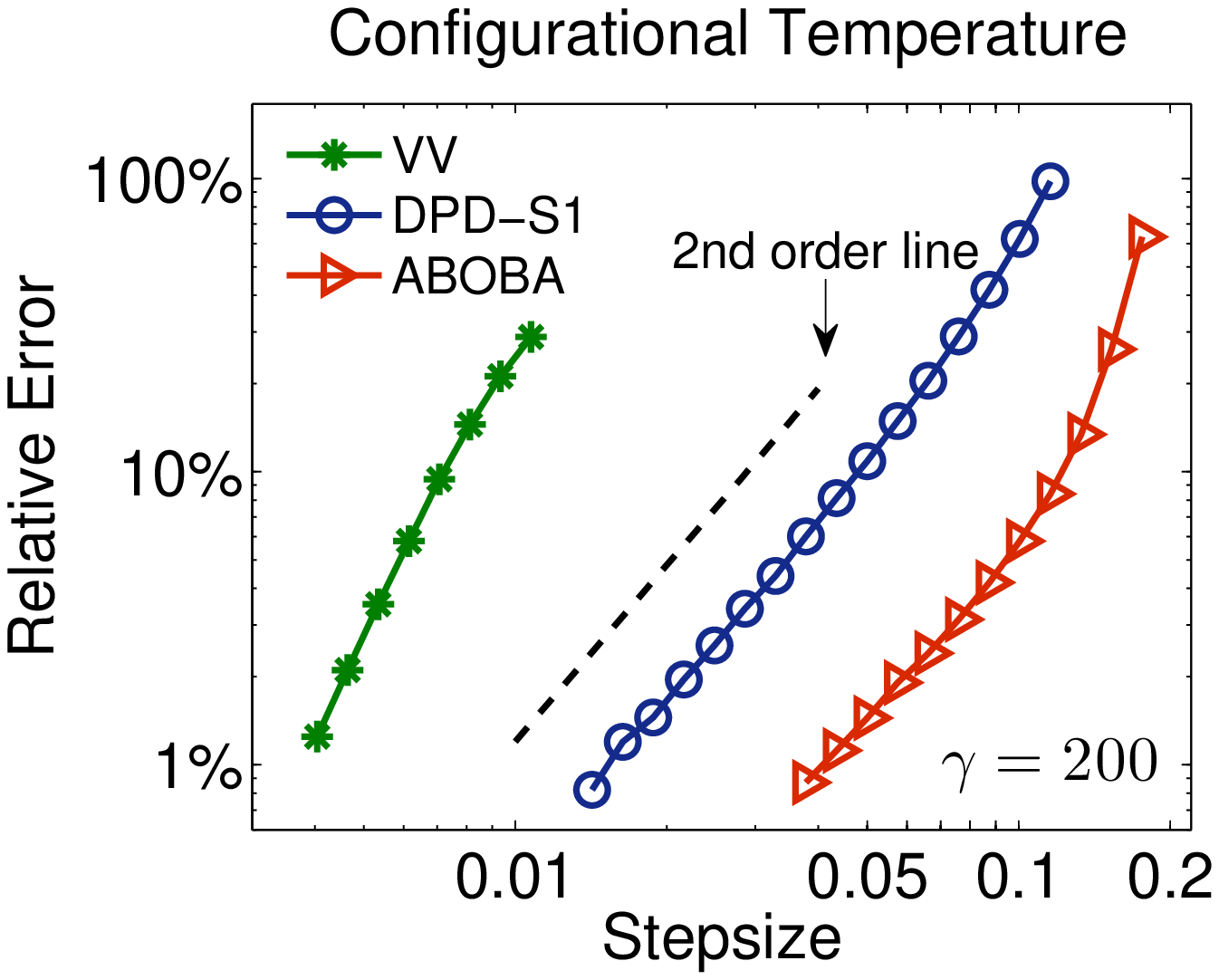}
\includegraphics[scale=0.4]{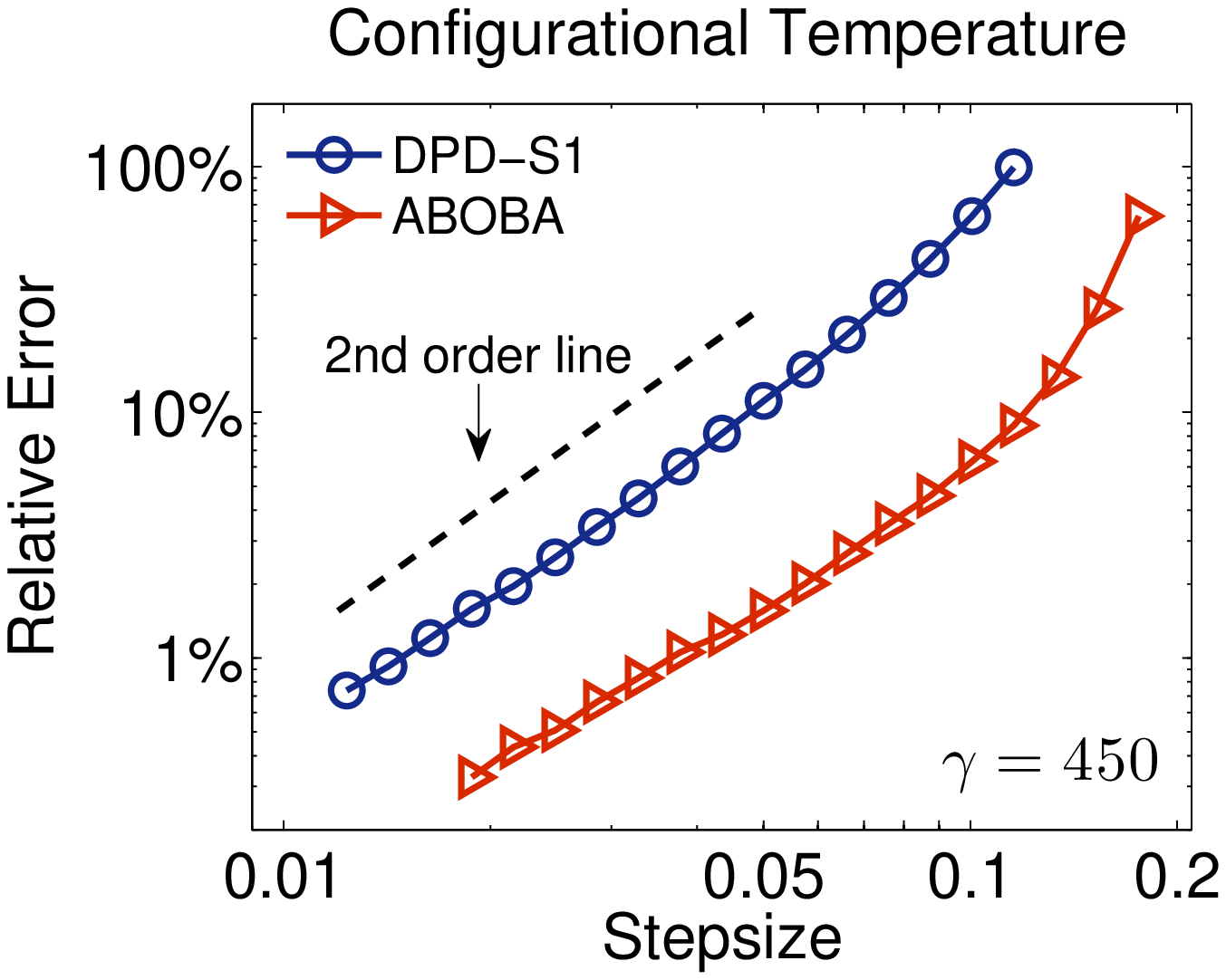}
\caption{\small Comparisons of the relative error in the computed configurational temperature~\eqref{eq:Config_Temp} against stepsize by using various numerical methods of the DPD system with friction coefficients of $\gamma=200$ (left) and $\gamma=450$ (right). The format of the plots is the same as in Figure~\ref{fig:DPD_CT_Comp_g4d5+g40d5}. Note that the VV method on the left panel became unstable at slightly over $\Delta t=0.01$.}
\label{fig:DPD_CT_Comp_g200+g450}
\end{figure}

The configurational temperature control for a variety of methods described in Section~\ref{sec:Numerical_Methods} was compared in Figure~\ref{fig:DPD_CT_Comp_g4d5+g40d5}. According to the dashed order line on both panels, we can see that all the methods tested exhibit second order convergence to the invariant measure, although only first order was expected for the VV method. Note that in some cases the errors appear to grow more rapidly when approaching their respective stability thresholds---notably the VV method on the right panel became unstable at slightly over $\Delta t=0.05$.

More specifically, in the friction limit of $\gamma=4.5$ on the left panel of Figure~\ref{fig:DPD_CT_Comp_g4d5+g40d5}, the DPD-S1 method is only marginally more accurate than the VV method with a fixed stepsize. Remarkably, the newly proposed ABOBA method achieves about one order of magnitude improvement over both methods in terms of the numerical accuracy with a fixed stepsize. The substantial improvement of ABOBA is quantified by measuring the numerical efficiency with a fixed level of accuracy in Table~\ref{table:efficiency_g4d5}. Note that the CPU times per step of all the three methods are similar to each other (particularly between DPD-S1 and ABOBA), indicating the methods have comparable complexities. (It is worth mentioning that, according to Table 2 in~\cite{Leimkuhler2015} via the CPU time of the DPD-S1 method, the complexity of ABOBA is also very similar to that of other DPD integrators, including DPD-Trotter~\cite{DeFabritiis2006,Serrano2006} and the Peters thermostat~\cite{Peters2004}.) However, with a fixed level of accuracy, the ABOBA method is able to use a much larger stepsize, thereby achieving a more than 104\% improvement over the benchmark VV method, whose numerical efficiency is very similar to that of DPD-S1. Note also that a 100\% improvement in the numerical efficiency effectively doubles the performance.

\begin{table}[t]
\centering
  \begin{tabular}{cccc}
  \hline
  \textbf{Method} & \textbf{Critical stepsize} & \textbf{CPU time} & \textbf{Scaled efficiency} \\ \hline
  VV     & 0.007 & 6.829 & 100.0\% \\
  DPD-S1 & 0.044 & 7.852 & 546.7\% \\
  ABOBA  & 0.116 & 7.742 & 1461.7\% \\ \hline
  \end{tabular}
\caption[Table caption text]{\small Comparisons of the ``numerical efficiency'' of various numerical methods of the DPD system with a very large friction coefficient of $\gamma=200$, corresponding to the left panel of Figure~\ref{fig:DPD_CT_Comp_g200+g450}. The format of the table is the same as in Table~\ref{table:efficiency_g4d5}.}
\label{table:efficiency_g200}
\end{table}

\begin{table}[t]
\centering
  \begin{tabular}{cccc}
  \hline
  \textbf{Method} & \textbf{Critical stepsize} & \textbf{CPU time} & \textbf{Scaled efficiency} \\ \hline
  DPD-S1 & 0.044 & 7.852 & 100.0\% \\
  ABOBA  & 0.116 & 7.742 & 267.4\% \\ \hline
  \end{tabular}
\caption[Table caption text]{\small Comparisons of the ``numerical efficiency'' of the DPD-S1 and ABOBA methods of the DPD system with an extremely large friction coefficient of $\gamma=450$, corresponding to the right panel of Figure~\ref{fig:DPD_CT_Comp_g200+g450}. The format of the table is the same as in Table~\ref{table:efficiency_g4d5} except the numerical efficiency of ABOBA was scaled to that of DPD-S1.}
\label{table:efficiency_g450}
\end{table}

The performance of various methods is largely similar in the large friction limit of $\gamma=40.5$ on the right panel of Figure~\ref{fig:DPD_CT_Comp_g4d5+g40d5} except the VV method appears to be slightly more accurate than DPD-S1 with a relatively small stepsize until it starts to blow up just over $\Delta t=0.03$. The ABOBA method again clearly outperforms both methods by at least one order of magnitude---Table~\ref{table:efficiency_g40d5} confirms that ABOBA achieves an almost 170\% improvement in the numerical efficiency over the benchmark VV method, while DPD-S1 only slightly outperforms VV.

We also investigate the performance of various methods with larger friction coefficients in Figure~\ref{fig:DPD_CT_Comp_g200+g450} where second order convergence to the invariant measure was again observed for all the methods tested except the VV method on the left panel, corresponding to the very large friction limit of $\gamma=200$. As in the case of $\gamma=40.5$ in Figure~\ref{fig:DPD_CT_Comp_g4d5+g40d5}, it is possible that the error appears to grow more rapidly as the VV method approaches its stability threshold just over $\Delta t=0.01$.

In the very large friction limit of $\gamma=200$ on the left panel of Figure~\ref{fig:DPD_CT_Comp_g200+g450}, with a fixed level of accuracy it can be easily seen that both DPD-S1 and ABOBA can use substantially larger stepsizes than that of VV. To be more precise according to Table~\ref{table:efficiency_g200}, a stepsize over six times as large as that of VV can be used for DPD-S1, contributing to a more than 440\% improvement in the numerical efficiency, while a stepsize over 16 times as large can be used for ABOBA, leading to a remarkable more than 1360\% enhancement. Due to its extremely poor performance in the very large friction limit of $\gamma=200$, in what follows VV will not be included for comparisons in that limit and beyond. In the extremely large friction limit of $\gamma=450$ on the right panel, ABOBA again comfortably outperforms DPD-S1, with a more than 160\% improvement according to Table~\ref{table:efficiency_g450}.

\subsubsection{Radial distribution function}
\label{subsubsec:RDF}

\begin{figure}[tb]
\centering
\includegraphics[scale=0.4]{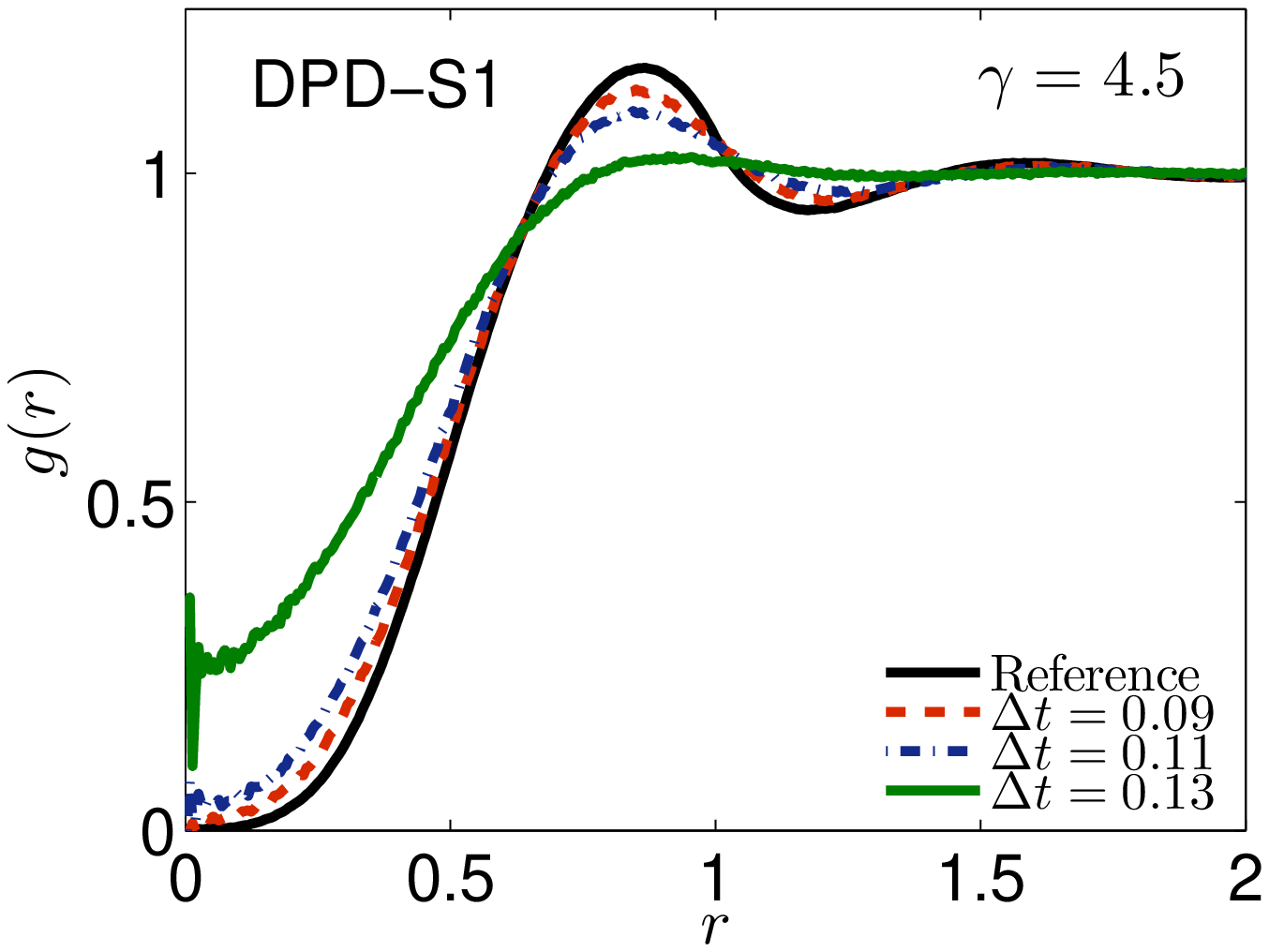}
\includegraphics[scale=0.4]{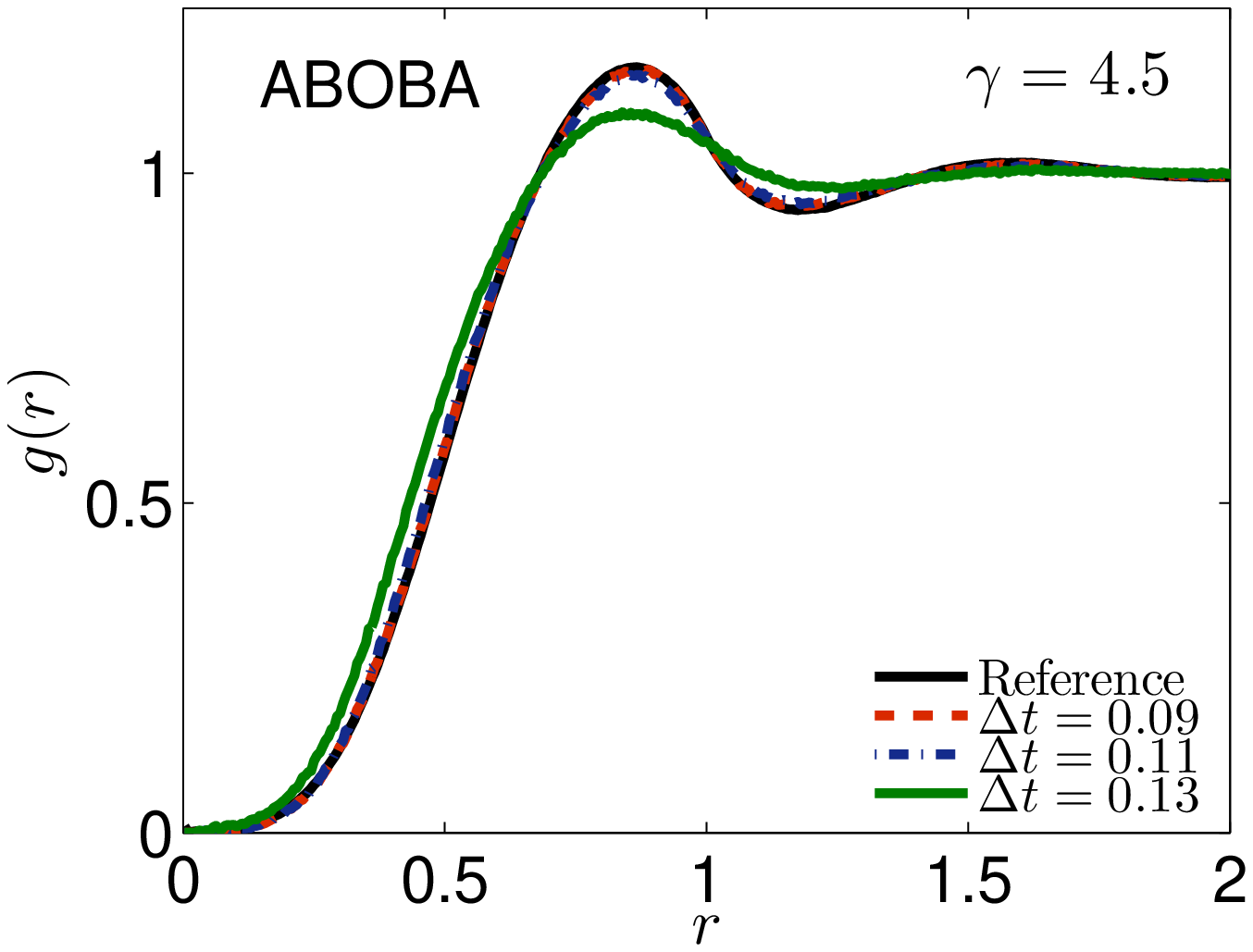}
\caption{\small (Color.) Comparisons of the radial distribution function (RDF), $g(r)$, obtained from the DPD-S1 method (left) and the ABOBA method (right) of the DPD system with a friction coefficient of $\gamma=4.5$. The solid black line is the reference solution obtained by using the DPD-S1 method with a very small stepsize of $\Delta t=0.001$, while the colored lines correspond to different stepsizes as indicated.}
\label{fig:DPD_RDF_Comp_g4d5}
\end{figure}

%
%
\begin{figure}[tb]
\centering
\includegraphics[scale=0.4]{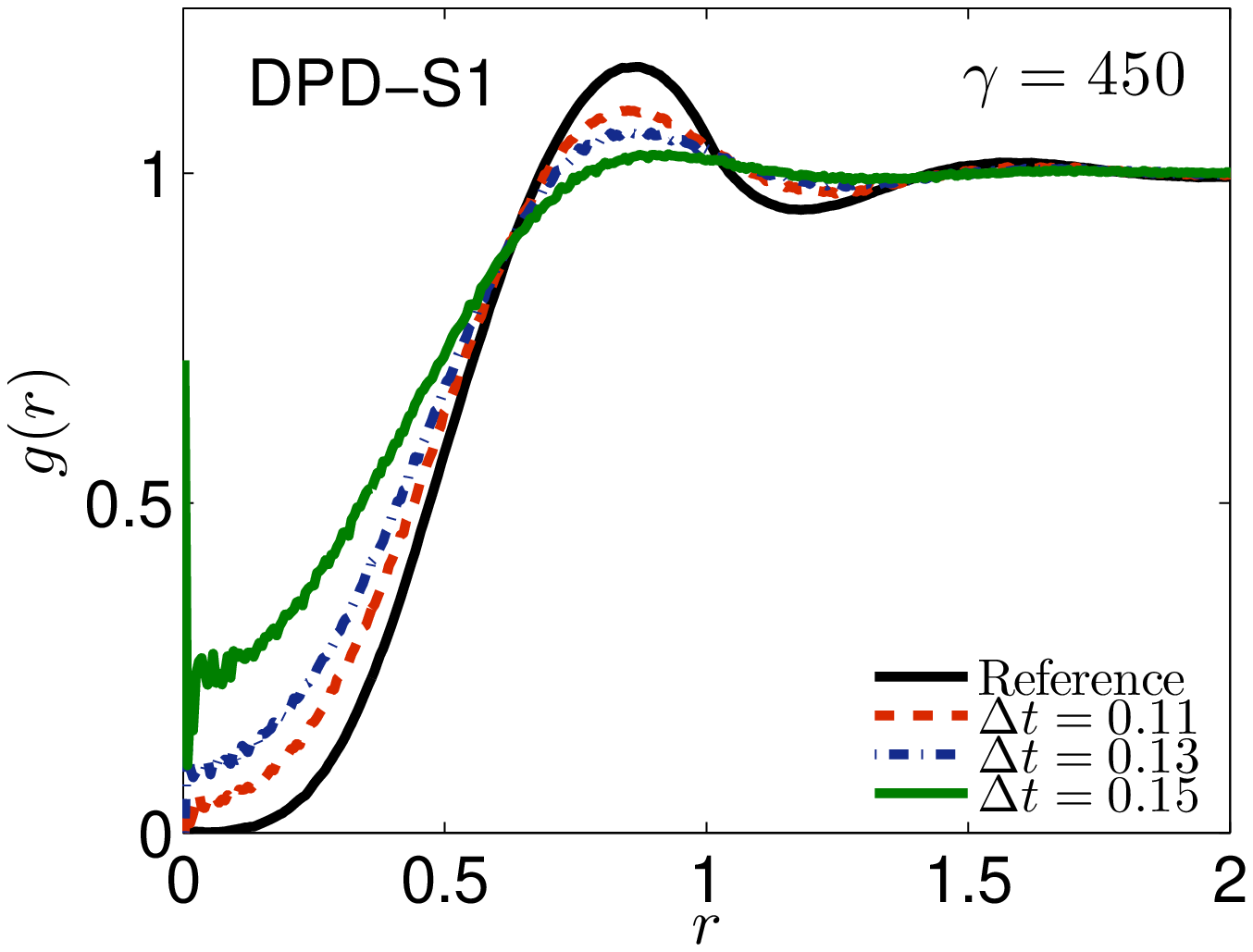}
\includegraphics[scale=0.4]{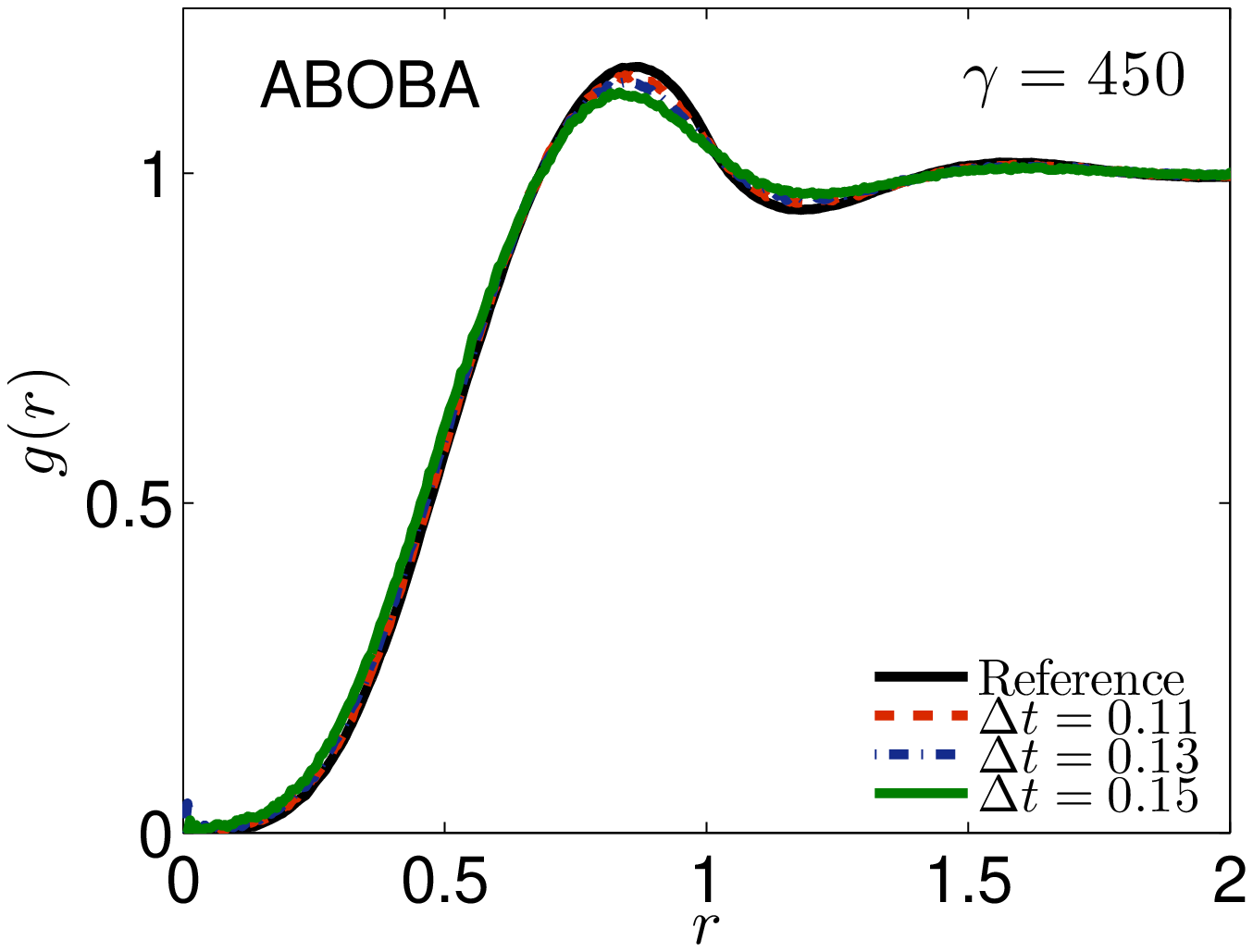}
\caption{\small (Color.) Comparisons of the radial distribution function (RDF), $g(r)$, obtained from the DPD-S1 method (left) and the ABOBA method (right) of the DPD system with an extremely large friction coefficient of $\gamma=450$. The format of the plots is the same as in Figure~\ref{fig:DPD_RDF_Comp_g4d5}.}
\label{fig:DPD_RDF_Comp_g450}
\end{figure}

Figure~\ref{fig:DPD_RDF_Comp_g4d5} compares the RDF that characterizes the structure of the DPD system with a friction coefficient of $\gamma=4.5$. The RDF of DPD-S1 appears to show pronounced artifacts at around $\Delta t=0.09$ and to be heavily destroyed at around $\Delta t=0.13$. (The RDF of the VV method in this case is very similar to that of DPD-S1 and thus we present only the results of DPD-S1 in Figure~\ref{fig:DPD_RDF_Comp_g4d5}.) Consistent with our findings on the left panel of Figure~\ref{fig:DPD_CT_Comp_g4d5+g40d5}, larger stepsizes can be used for ABOBA without compromising the control of the structure of the system---with a stepsize as large as $\Delta t=0.11$, the RDF on the right panel of Figure~\ref{fig:DPD_RDF_Comp_g4d5} is almost indistinguishable from the reference solution, while only small deviations were observed in the RDF with a stepsize of $\Delta t=0.13$.


Since the RDF controls of both DPD-S1 and ABOBA methods were very similar with a wide range of large friction coefficients (e.g., $\gamma=40.5$, $\gamma=200$, and $\gamma=450$), only the results with the largest friction coefficient of $\gamma=450$ were included in Figure~\ref{fig:DPD_RDF_Comp_g450}. Again consistent with our findings on the right panel of Figure~\ref{fig:DPD_CT_Comp_g200+g450}, larger stepsizes can be used for ABOBA while maintaining good control of the structure of the system---while, with a stepsize of $\Delta t=0.15$, the RDF of ABOBA shows rather small deviations from the reference solution, it is heavily destroyed for DPD-S1. It is worth mentioning that in the large friction limit of $\gamma=40.5$ (results not shown), the performance of the methods also largely align with our observations on the right panel of Figure~\ref{fig:DPD_CT_Comp_g4d5+g40d5}, particularly, the VV method became unstable at slightly over $\Delta t=0.05$, at which point the RDF visibly deviates from the reference solution.

As expected in Section~\ref{subsec:ABOBA} and observed numerically, the computational cost (i.e., CPU time per step) of the BAOAB method in the DPD context is almost twice as much as that of the ABOBA method. Therefore, unless at least twice as large a stepsize can be used in BAOAB as that of ABOBA while maintaining the same level of accuracy in physical quantities of interest, the numerical efficiency of the former will be outperformed by the latter. We observed that the performance of BAOAB in terms of RDF is very similar to that of ABOBA. Although BAOAB is able to outperform ABOBA in the numerical efficiency based on the computed configurational temperature with friction coefficients of $\gamma=200$ and $\gamma=450$ if the level of accuracy required is low enough (e.g., less than 1\%), overall we conclude that the ABOBA method should be preferred. It is also worth mentioning that a superconvergence (i.e., a fourth order convergence to the invariant measure) result observed only in the BAOAB method in the large friction limit in Langevin dynamics~\cite{Leimkuhler2013,Leimkuhler2013c} (see another example in adaptive Langevin dynamics~\cite{Leimkuhler2015a}) was not observed in the DPD context. In order to ensure the BAOAB method is symmetric, one may modify the method to be an BAOOAB-like method where notably the order of the interacting pairs in either of the O parts needs to be reversed. However, in our numerical experiments we observed that the performance of the BAOOAB method was very similar to that of the BAOAB method---the superconvergence property was not observed in the BAOOAB method either. It is also worth pointing out that solely replacing the BBK method by the analytical solutions in the OU process of the DPD-S1 method has little impact on its performance, which indicates that it is the ordering of the pieces in the ABOBA method that makes more of the difference; nevertheless, the analytical solutions should be preferred.

\subsubsection{Shear viscosity}
\label{subsubsec:Shear_Viscosity}

\begin{figure}[t]
\centering
\includegraphics[scale=0.4]{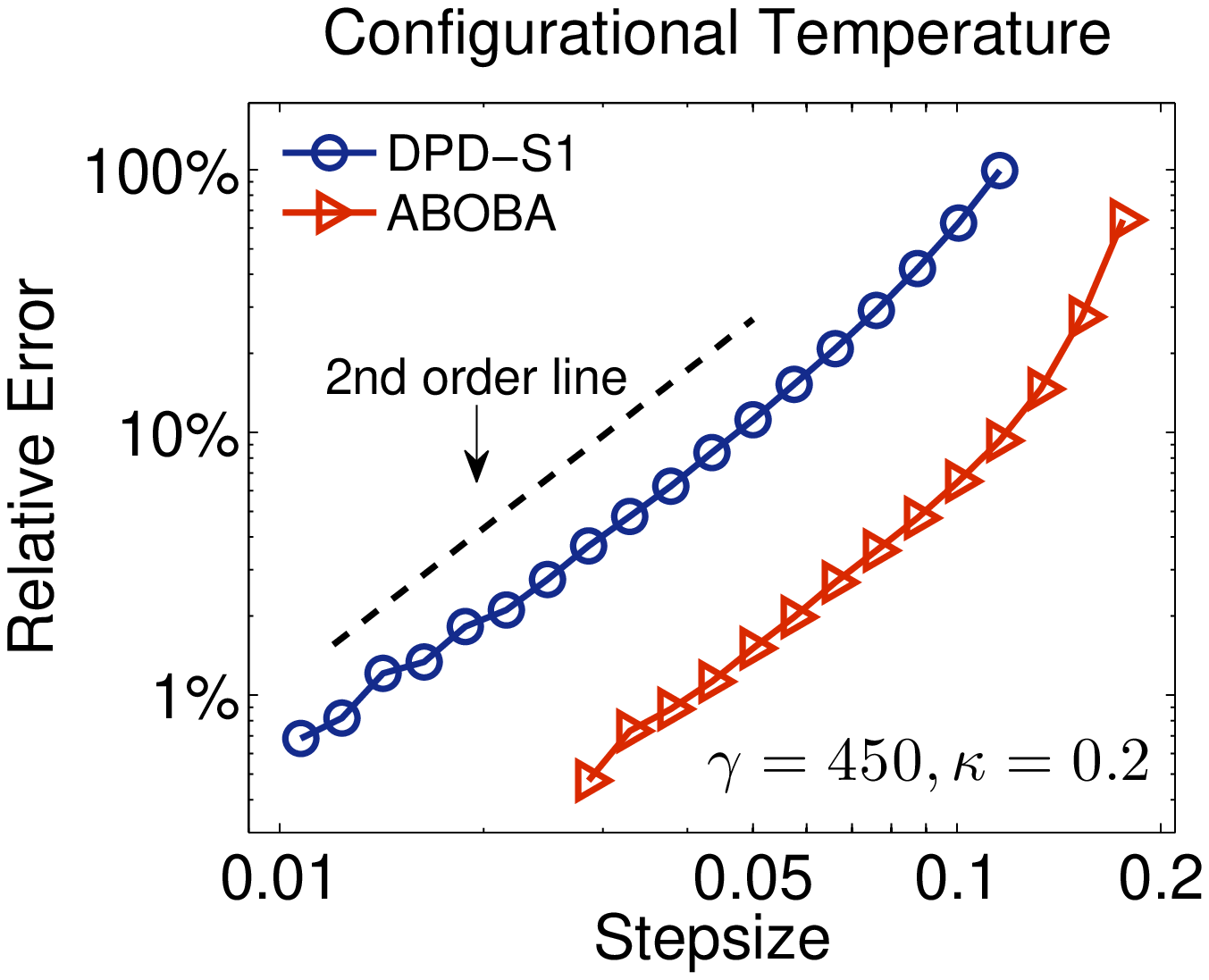}
\includegraphics[scale=0.4]{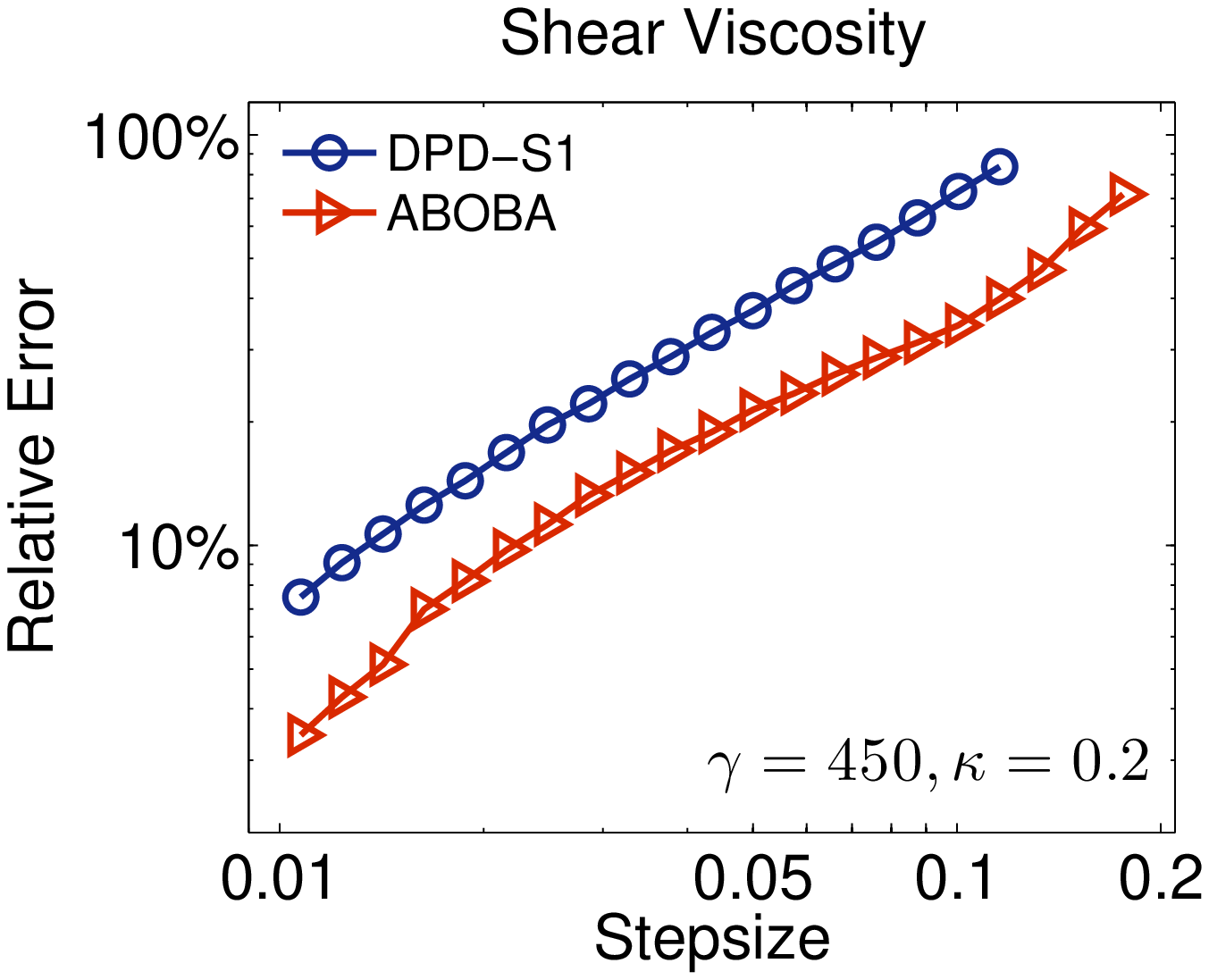}
\caption{\small Comparisons of the relative error in the computed configurational temperature (left) and shear viscosity (right) against stepsize by using both DPD-S1 and ABOBA methods of the DPD system with a friction coefficient of $\gamma=450$ and a shear rate of $\kappa=0.2$ generated by the Lees--Edwards boundary conditions. The format of the plots is the same as in Figure~\ref{fig:DPD_CT_Comp_g4d5+g40d5}. Note that in both cases the reference values were obtained by using the DPD-S1 method with a very small stepsize of $\Delta t=0.001$. }
\label{fig:DPD_CT+Vis_Comp_g450_kd2}
\end{figure}

Since we are more interested in a fluid-like Schmidt number when measuring transport coefficients, we restrict our attention to the extremely large friction limit of $\gamma=450$ in this section. The Lees--Edwards boundary conditions were applied in order to generate a simple and steady shear flow with a shear rate of $\kappa=0.2$. Note that since we are interested in the Newtonian regime, the shear rate imposed should be sufficiently small as suggested in~\cite{Pagonabarraga1998}.

The configurational temperature control for both DPD-S1 and ABOBA methods was assessed on the left panel of Figure~\ref{fig:DPD_CT+Vis_Comp_g450_kd2}. The performance of both methods is very similar to the equilibrium case on the right panel of Figure~\ref{fig:DPD_CT_Comp_g200+g450}. With a 10\% relative error in the configurational temperature, over twice as large a stepsize can be used in ABOBA as that of DPD-S1, resulting in a more than 135\% improvement in the numerical efficiency according to Table~\ref{table:efficiency_g450_kd2_CT}.

The right panel of Figure~\ref{fig:DPD_CT+Vis_Comp_g450_kd2} compares the performance of both methods in terms of the shear viscosity control. The ABOBA method again clearly outperforms the DPD-S1 method in terms of the numerical accuracy with a fixed stepsize. With a 10\% relative error in the shear viscosity, almost twice as large a stepsize can be used in ABOBA as that of DPD-S1, leading to an impressive almost 86\% improvement in the numerical efficiency according to Table~\ref{table:efficiency_g450_kd2_Vis}.

\section{Conclusions}
\label{sec:Conclusions}

\begin{table}[t]
\centering
  \begin{tabular}{cccc}
  \hline
  \textbf{Method} & \textbf{Critical stepsize} & \textbf{CPU time} & \textbf{Scaled efficiency} \\ \hline
  DPD-S1 & 0.050 & 7.852 & 100.0\% \\
  ABOBA  & 0.116 & 7.742 & 235.3\% \\ \hline
  \end{tabular}
\caption[Table caption text]{\small Comparisons of the ``numerical efficiency'' of the DPD-S1 and ABOBA methods of the DPD system with an extremely large friction coefficient of $\gamma=450$ and a shear rate of $\kappa=0.2$, corresponding to the left panel of Figure~\ref{fig:DPD_CT+Vis_Comp_g450_kd2}. The format of the table is the same as in Table~\ref{table:efficiency_g450}. }
\label{table:efficiency_g450_kd2_CT}
\end{table}

\begin{table}[t]
\centering
  \begin{tabular}{cccc}
  \hline
  \textbf{Method} & \textbf{Critical stepsize} & \textbf{CPU time} & \textbf{Scaled efficiency} \\ \hline
  DPD-S1 & 0.012 & 7.852 & 100.0\% \\
  ABOBA  & 0.022 & 7.742 & 185.9\% \\ \hline
  \end{tabular}
\caption[Table caption text]{\small Comparisons of the ``numerical efficiency'' of the DPD-S1 and ABOBA methods of the DPD system with an extremely large friction coefficient of $\gamma=450$ and a shear rate of $\kappa=0.2$, corresponding to the right panel of Figure~\ref{fig:DPD_CT+Vis_Comp_g450_kd2}. The format of the table is the same as in Table~\ref{table:efficiency_g450} except the physical quantity of interest is the shear viscosity. }
\label{table:efficiency_g450_kd2_Vis}
\end{table}

We have proposed one promising splitting method of the DPD system, namely the ABOBA method, which relies on solving each of the splitting parts exactly. We have also demonstrated and numerically verified the second order convergence to the invariant measure for ABOBA. Although only first order convergence was expected for the VV method, second order convergence was observed in our numerical experiments. While the VV method performs comparably with the DPD-S1 method with a friction coefficient of $\gamma=4.5$, a standard choice in DPD simulations, it has been demonstrated that it is not suitable for (much) larger friction coefficients required for a fluid-like Schmidt number. Remarkably, the newly proposed ABOBA method substantially outperforms both VV and DPD-S1 methods with a wide range of the friction coefficients in all the physical quantities tested.

To be more precise, in terms of the configurational temperature control, ABOBA is at least one order of magnitude more accurate than the benchmark VV method with a fixed stepsize in each of the cases in Figures~\ref{fig:DPD_CT_Comp_g4d5+g40d5}--\ref{fig:DPD_CT_Comp_g200+g450}. Moreover, VV became unstable easily (slightly over $\Delta t=0.05$ with $\gamma=40.5$ and $\Delta t=0.01$ with $\gamma=200$) as the friction coefficient was increased, which indicates that vanishingly small stepsizes have to be used in the extremely large friction limit of $\gamma=450$, effectively ruling it out especially for large-scale DPD simulations with a fluid-like Schmidt number. Furthermore, ABOBA constantly outperforms the most recommended DPD-S1 method, achieving a more than 100\% improvement (i.e., at least doubling the performance) in each of the cases in Tables~\ref{table:efficiency_g4d5}--\ref{table:efficiency_g450_kd2_CT}. As previously mentioned, good control of the configurational temperature appears to lead to good performance in other physical quantities. In all the comparisons of the radial distribution function, ABOBA again can use larger stepsizes than that of DPD-S1 in order to maintain the same level of accuracy. The performance of the ABOBA method is equally impressive in our nonequilibrium simulations driven by Lees--Edwards boundary conditions~\cite{Lees1972}; it achieves an almost 86\% improvement over DPD-S1 in the numerical efficiency based on the computed shear viscosity, an important transport coefficient in nonequilibrium molecular dynamics. This again illustrates the importance of optimal design of numerical methods.

Although the VV method has been substantially outperformed by the ABOBA method, the parallelization of the latter (along with DPD-S1), due to the fact that the interacting pairs in the OU process are solved successively, is not as straightforward as that of the former. However, a similar task of parallelizing the DPD-S1 method has been addressed by Larentzos et al.~\cite{Larentzos2014}---the procedures can be easily adopted for parallelizing ABOBA.

Inspired by recent developments in adaptive thermostats~\cite{Jones2011,Leimkuhler2015a,Shang2015} (their theoretical foundations have recently been taken up in~\cite{Leimkuhler2019a}), the so-called pairwise adaptive Langevin (PAdL) thermostat has been proposed in~\cite{Leimkuhler2016a}. PAdL is able to correct for thermodynamic observables while mimicking the dynamical properties of DPD and thus can be viewed as ``adaptive DPD''. It has been demonstrated in~\cite{Leimkuhler2016a} that PAdL can also outperform popular numerical methods (including DPD-S1) for DPD in both equilibrium and nonequilibrium settings. Since PAdL is not based on the standard DPD formulation, it has not been included for comparisons in the current article. However, it will be interesting to compare its performance with the outstanding ABOBA method proposed in the current article. We leave further exploration of this direction for future work. As Shardlow's splitting method, the newly proposed ABOBA method can also be applied in DPD with various fixed conditions~\cite{Lisal2011}; it is worth assessing the performance of the ABOBA method in those settings. It is also worth investigating the performance of various methods in the case of DPD with constraints~\cite{Gavrilov2020}.

\section*{Acknowledgements}

The author thanks Richard Anderson, Michael Seaton, Patrick Warren, and the anonymous referees for their valuable suggestions and comments.

\bibliographystyle{is-abbrv}

\bibliography{refs}

\end{document}